\documentclass[preprint,review,10pt]{elsarticle}

\usepackage{amssymb}
\usepackage{amsthm}
\usepackage{amsmath,amsfonts}
\usepackage{algorithmic}
\usepackage{algorithm}
\usepackage{array}
\usepackage[caption=false,font=normalsize,labelfont=sf,textfont=sf]{subfig}
\usepackage{textcomp}
\usepackage{stfloats}
\usepackage{url}
\usepackage{verbatim}
\usepackage{graphicx}
\usepackage{natbib}
\usepackage{lineno}
\usepackage{color}
\usepackage[dvipsnames]{xcolor}
\usepackage{bm}
\newcommand\ie{\textit{i.e.}}
\newcommand\eg{\textit{e.g.}}

\journal{Pattern Recognition}

\begin{document}

\begin{frontmatter}



\title{Pattern based learning and optimisation through pricing for bin packing problem}

\author[unnc]{Huayan Zhang}
\author[unnc]{Ruibin Bai\corref{cor1}}
\author[msra]{Tie-Yan Liu}
\author[unnc]{Jiawei Li}
\author[unnc]{Bingchen Lin}
\author[unnc]{Jianfeng Ren}

\cortext[cor1]{Corresponding Author: ruibin.bai@nottingham.edu.cn}

\affiliation[unnc]{organization={School of Computer Science},%Department and Organization
            addressline={University of Nottingham Ningbo China}, 
            city={Ningbo},
            state={Zhejiang},
            country={China}}
\affiliation[msra]{organization={Microsoft Research AI4Science},
            addressline={Microsoft Research Asia},
            state={Beijing},
            country={China}}
\begin{abstract}

As a popular form of knowledge and experience, patterns and their identification have been critical tasks in most data mining applications. However, as far as we are aware, no study has systematically examined the dynamics of pattern values and their reuse under varying conditions. We argue that, when problem conditions such as the distributions of random variables change, the patterns that performed well in previous circumstances may become less effective and adoption of these patterns would result in sub-optimal solutions.
In response, we make a connection between data mining and the duality theory in operations research and propose a novel scheme to efficiently identify patterns and dynamically quantify their values for each specific condition. Our method quantifies the value of patterns based on their ability to satisfy stochastic constraints and their effects on the objective value, allowing high-quality patterns and their combinations to be detected. We use the online bin packing problem to evaluate the effectiveness of the proposed scheme and illustrate the online packing procedure with the guidance of patterns that address the inherent uncertainty of problem. Results show that the proposed algorithm significantly outperforms the state of the art methods. We also analysed in detail the distinctive features of the proposed methods that lead to the performance improvement and the special cases where our method can be further improved.  
\end{abstract}


\begin{highlights}
    \item Generalise pattern as a form of knowledge but with changing values under different conditions. 
    \item A novel mechanism is proposed to accurately quantify the prices of patterns under known distributions. 
    \item We extended it further with an adaptive predictive-reactive framework for unknown distributions.
    \item Our method significantly outperforms the existing online bin packing methods.
    % \item Generalise the concept of pattern and associated plan to online bin packing problems.
    % \item Identifying efficient patterns by dynamic pricing mechanism.
    % \item Utilises measurement of uncertainty for online decision making to eliminate imperfect planning.
    % \item A learning-based mechanism to automatically detect item distribution for online planning.
    % \item Outperform mostly used heuristics and DRL benchmark model without full knowledge of problem instance.
    % \item Outperforms recent SOTA online algorithms without full knowledge of problem instance.

\end{highlights}


\end{frontmatter}

\linenumbers

% \section{}
% \label{}
\section{Introduction}
Combinatorial optimisation problems (COP) have extensive applications in various industrial fields \cite{korte2011combinatorial}. % \rjf{One reference to a survey paper? Or a few references here.} 
However, due to the NP-hardness nature of such problems, finding optimal solutions becomes extremely challenging given limited computational power, particularly for large-scale instances.
% \rjf{Add references to support this statement.} 
This challenge escalates further when uncertainties are considered %\rjf{refs.}, 
which hinders us from deriving the practical solutions.

Existing approaches for addressing these types of problems can be broadly categorised into two main groups~\citep{bai-analytics-2023}: analytical model-driven methods, which are often exemplified by analytical and mathematical models ~\citep{gupta2020interior}% \red{analysis techniques}
 % \rjf{This term is too general. Can you be more specific?}
; and data-driven methods such as genetic programming and reinforcement learning~\citep{chen2024deep}%~\citep{BENGIO2021405, ULMER2020100008}.  
The former primarily concentrates on the analytical properties of problem models, but it may encounter challenges in terms of robustness when confronted with uncertainties in the input data~\citep{bai-analytics-2023}. In contrast, data-driven methods typically approach combinatorial problems as online optimisation problems %\rjf{ref}. 
They often address the problem sequentially, employing policies or rules that account for the realisation of random variables and the states of the partial solution at each decision point.
One of the major limitations of data-driven methods is their inability to effectively exploit the core structure and properties of the problem \cite{bai-analytics-2023}. Specifically, existing data-driven approaches~\citep{Lu2020A} often prioritise optimisation objectives while neglecting the intricate inter-dependencies among decision variables (represented as constraints), and their cumulative influence on the overall objective. 
 % \rjf{One or two sentences summarize the current challenges and describe their relations to pattern.} 

Patterns are one of the most powerful and effective problem-solving tactics in computer vision~\citep{GOMEZ20149, CHEN2022108769,LIN2024110143} and time-series data analysis~\citep{BAI202046,LI2021107711, BREITENBACH2023109355}. 
%It has been widely used in solving combinatorial optimisation problems like bin packing~\citep{angelopoulos2022online,sym14071301,BORTFELDT2020545} and vehicle routing~\citep{Lu2020A,ARNOLD2021107957}. Patterns can be in diverse forms across different fields. 
However, very limited number of studies have used patterns as a problem-solving strategy for combinatorial optimisation. A pioneer work was made in 2021 by a group of academics from MIT and AI experts from Amazon in their Last-mile Routing Research Challenge\footnote{see details from https://routingchallenge.mit.edu}, whose primary objective is to search for ``high-valued'' vehicle route patterns that take into account not just the route lengths/costs, but also the tacit knowledge linked to safety, robustness and sustainability. However, in this competition, the identification of the route patterns still relies on manually labelling the qualify of large number of routes which are expensive and not transferable to other cases. More importantly, robustness of these labels becomes questionable when confronted with uncertainties. 

As a problem-solving strategy, pattern-based method has several advantages in solving COPs. 
% \rjf{Suggest to come back and restrict the discussions to COPs here.} 
First, patterns are interpretable, editable and reusable %\rjf{refs}.
Therefore, the solutions built from patterns can be better comprehended and easier adapted to practical applications. 
Second, patterns allow complex (including non-linear) constraints to be modelled implicitly. For example, in the aforementioned Amazon's last-mile routing problem, drivers' duty obligations and preferences are embedded in vehicle route patterns. This enables us to build a pattern-based linear problem model because nonlinear constraints are handled within the pre-generated patterns.% ~\citep{Scheithauer2018}.  
 % \rjf{This example is not intuitive enough for readers to understand the second advantage. By the way, the second advantage itself is not intuitive enough.} 
Finally, patterns can be considered as a form of knowledge or experience-originated rules that can be analysed and migrated to new problems of similar structure, promoting transfer learning and knowledge reuse, albeit under the guidance of their dynamic values. 
%that vary in distribution~\citep{burke_providingmemorymechanism_2010}. 
%  and polished to build better understanding about the problems and their solutions~\citep{burke_providingmemorymechanism_2010}. 
% \rjf{What is the difference from advantage 1: interpretable? Maybe here we emphasise reusable and remove reusable in advantage 1?}. \rjf{One or two sentence summarise in which conditions patterns may work well and the benefits of using patterns.} 

For many combinatorial optimisation problems with uncertainties, it is challenging to derive good patterns that work well across different scenarios. 
%the identification of good patterns alone is not sufficient to solve the problems. 
%When the problem input data carry certain degree of uncertainties, 
Patterns that are considered ``good" in some contexts may become of poorer quality due to uncertainties in different problem-solving scenarios. The underlying cause is that, under uncertainties, although the problem structure (in terms of the objective function and the constraint structures) remains unchanged from instance to instance, the inter-dependencies between the decision variables may have changed significantly, leading to performance drop when utilising some of the patterns. It is therefore critical to find a way to accurately quantify the value of patterns under different problem-solving conditions so that the most suitable patterns can be derived adaptively for different stochastic scenarios.

%We shall illustrate this in the next section through a numerical example. 

%In this research, 
%This research focuses the online 1D bin packing problem, as it is one of the most studied online COP, where decision must be made immediately once an item arrives~\cite{10.1007/978-3-642-58409-1_2}. The future problem state are either partially visible~\cite{9982095} or invisible~\cite{8956393} to the decision maker. As the consequence, solutions made by online decision process are usually short-sighted and semi-optimal. A natural intuition is to extract patterns and associated combinations based on historical observation, as a reference for allocating future items \cite{angelopoulos2022online}. However, due to the problem information is incomplete, future problem state might be vary from past observation. Moreover, sometimes the problem is stochastic in terms of quantity of each item types. These issues may affect how valuable a pattern is in online decision process.
To address the aforementioned challenges, in this paper, we propose a novel scheme that can systematically generate high-valued patterns for each perceived stochastic scenario and then optimise their reuse in the near-optimal way. Our method is built on the concept of duality and shadow prices in linear programming. The effectiveness of the proposed method is examined on 1D online bin packing problem, one of the most intensively studied COPs with many practical applications. 
The shadow price depicts the marginal impact of constraints on the objective function. In the one-dimensional bin packing problem, this reflects the change in the optimal solution when the number of a certain item in the sequence changes. By calculating shadow prices, we can dynamically determine the importance of items for different distributions. Mathematically, shadow price is usually obtained by solving the dual problem of the original problem. 
% The value of potential patterns can be estimated by shadow prices 
Guided by the shadow prices, patterns with potential to improve the objectives are repeatedly generated through solving pricing sub-problems. This generation method is also referred to as \textit{column generation}, in which a pattern corresponds a column in the left-hand matrix of the pattern-based linear programming formulation. The column generation process is stopped when no new pattern can be found to improve the objective value, leading to the optimal solution defined by a combination of adopted patterns and their frequency of use.  
For online problems, the optimal pattern combination is generated by the above method based on the latest forecast.
The patterns combination is used as a packing plan to guide online packing procedure. Due to the inherent uncertainty and imprecise forecasts, the online packing procedure must dynamically adjust packing plan by tracking the uncertainty during packing. 
% The whole online bin packing solving framework is shown in Figure~\ref{fig:framework}.
The proposed framework of solving online bin packing problem is named as Column Generation Plan-and-Pack (\texttt{CGPP}), shown in Figure~\ref{fig:framework}.
Our contributions can be summarised as follows: 
\begin{itemize}
    \item This paper introduces a novel scheme, namely \texttt{CGPP}, % \rjf{What algorithm? the name?} 
    to the field of learning-based online optimisation. By leveraging mathematical rigour of the duality theory and shadow prices in operations research, our method can discover high-quality reusable patterns while accurately quantifying their values for known uncertainty distributions, leading to near-optimal solutions and significant performance improvements compared to the current methods for online bin packing%\rjf{How to achieve this? The technique?} 
    
    \item We further extend the algorithm with an adaptive scheme for instances with unknown uncertainty distributions. Thanks to the improved ability of distribution forecasting of the scheme and its advanced packing strategies with imperfect plans, the resulting algorithm achieved outstanding results and significant improvement over existing methods. 

    \item The resulting high-level information and metrics (e.g. patterns and their use frequency and values) from our method provide deeper understanding and insights of the problem instances being addressed. Resulting solutions are constructed by fully understandable blocks in patterns, and hence can expect good acceptance by practitioners because of its obvious interpretabability. The method is generalisable to other real-life problems with similar structures. 
    
    % firstly using pricing technique, specifically column generation, to identify patterns in packing procedure. By applying pricing technique, we achieved significant advance comparing with other online decision algorithms like best-fit and reinforcement learning.
    
    %We introduced measurement of uncertainty in the online 1d bin packing problem, namely underestimate and overestimate. An adaptive packing procedure that adjust the plan when the incoming item is biased with prediction is applied to eliminate the gap between plan and reality.

\end{itemize}

\section{Preliminaries}
\subsection{Bin packing problem (BPP)}
The bin packing problem is formally defined as packaging a set of items of different sizes using the minimum number of boxes of the same capacity. In its basic offline version, The size of items is given before packaging. Let $B$ denote the capacity of the bins to be used and $T$ be the number of item types, with each item type $t$ having a size $s_t$ and quantity $q_t$. 
Let $y_j$ be a binary variable to indicate whether bin $j$ is used in a solution ($y_j=1$) or not ($y_j=0$). Let $x_{tj}$ be the number of times item type $t$ is packed in bin $j$. The problem can be formulated by

\begin{eqnarray}
\label{eqn:packing}
\mbox{minimise} \quad    & \sum_{j=1}^U y_{j} &\label{eqn:objective} 
\\
\mbox{subject to:} \quad & \sum_{j=1}^U x_{tj} = q_t &  \mbox{ for } t=1,\cdots,T  \label{eqn:demand} 
\\
& \sum_{t=1}^T s_t x_{tj} \leq B y_{j}      &  \mbox{ for } j=1,\cdots,U \label{eqn:constraint} 
\end{eqnarray}
where $U$ is the maximum number of possible bins that can be used. Bin packing problem is proven NP-Hard~\citep{martello_lowerboundsreduction_1990}. To improve the computational efficiency, heuristic and meta-heuristics are often used. The most well-known heuristics include Best Fit (BF), Minimum Bin Slack (MBS) and their variants~\citep{bai-simulated-2012}. 

\subsection{Online bin packing problem}
Although most research efforts on BPP have been focusing on its \textit{offline} version in which details of items to be packed are perfectly predictable in advance, many real-life packing problems appear to be \textit{online} because of dynamic realisation of items' specifications. More specifically, in model defined in Eq.~\eqref{eqn:packing}, the quantity of item type $t$, $q_t$, is often unknown but its proportion among all item types can be estimated. Items arrive sequentially over time and its information (\ie, the type of the current item) is only available after their arrivals. A solution method for online BPP must assign a bin to each randomly arrived items upon its arrival and this assignment cannot be subsequently altered. Therefore, best fit remains a good solution method for online BPP with a high competitive ratio (\eg, 1.7)~\citet{Scheithauer2018}, but MBS is not usable anymore because it relies on the full information of items to be packed. In this research, we aim to improve average performance for online BBP by solving it with a pattern based optimisation scheme guided by \textit{pricing}. The motivation and underlying ideas can be illustrated by using two simple BPPs given in Table~\ref{tab:bpp} with bin capacity $B=10$. 
\begin{table}[htbp]
    \centering
    \caption{Two simple 1D BPPs. $s_t$ are sizes of item types and bin capacity $B=10$.}

    \begin{tabular}{l|l}
    \hline
    case 1 & $\{s_t\}= \{5,4,4,3,2,2\}$ \\
    & best fit solution: $\{5,4\}, \{4,3,2\},\{2\}$\\
    & opt. solution: $\{5,3,2\}, \{4,4,2\}$\\
    \hline
    case 2 & $\{s_t\}= \{5,4,4,3,3,3,3,3,2,2,2,2,2,2\}$ \\
    & best fit solution: $\{5,4\},\{4,3,3\},$ \\ 
    & $\{3,3,3\},\{2,2,2,2,2\},\{2\}$ \\
    & opt. solution: $\{5,3,2\}, \{4,3,3\},$\\
    &$\{4,3,3\},\{2,2,2,2,2\}$\\
    \hline 
    \end{tabular} 
    \label{tab:bpp}

\end{table}

For both cases, best fit produces sub-optimal solutions. Like most heuristic methods, best fit is a typical objective-focused incremental method that aims to obtain the maximum possible benefits in terms of the objective defined in Eq.~\eqref{eqn:packing} at every step. However, although it partially addresses the constraint in Eq.~\eqref{eqn:constraint} by using a simple rule to seek the best possible packing, it fails to address the constraint in Eq.~\eqref{eqn:demand} completely because it does not proactively considers the quantities of each item. Indeed, this is one of the main problems of many existing learning-based approaches such as genetic programming (GP) and deep reinforcement learning (DRL), \eg, it is challenging to address multiple constraints while optimising the learning objective. 

The second challenge is that the packing patterns in the optimal solutions change significantly, as shown in Table \ref{tab:bpp}. Only one packing pattern in Case 1 is reused in Case 2, \ie, $\{5,3,2\}$. The optimal packing pattern $\{4,4,2\}$ in Case 1 disappeared completely in Case 2 and two new packing sets are now introduced for Case 2. Evidently, the previously good packing sets may not be good any more for the new instance while previously unpopular patterns may become more valuable. This imposes great challenges to algorithms that aim to exploit good patterns for solving the bin packing problems. 

The third challenge is that although the two instances are similar in terms of problem structure, \eg, bin size, item types, but the distribution of item types is different. In Case 2, more small items need to be packed. Therefore, the algorithm would need to not only look at how well each packing set performs in terms of the capacity waste, but also consider %\red{how efficient it packs items that appear more frequently}. 
how efficient it satisfies the demands of different item types in the incoming sequence. 
% \rjf{In another word, for online bin packing, the algorithm should not only consider how to optimally pack the existing items, it should also consider the distribution of items and predict the incoming items so that both the existing items and the possible incoming items could be optimally packed.} 

% \rjf{Check the correctness.}
In another word, for online bin packing, the algorithm should consider not only the existing (partial) packing state, but also the incoming item sequence so that packing results can be optimised on a longer-term scale.
% \myresponse{A little bit weird expression, we cannot pack existing items optimally, the decision is made immediately.}
% \rjf{For items or patterns that rarely appear, the algorithm should avoid to pack a bin following these patterns, as it is unlikely that the bin will be fully packed following these rare patterns.} \rjf{Check correctness.}
As the consequence, the algorithm is likely to have a better packing performance for long sequences, comparing with short-sighted methods.
% \myresponse{Weird too. Better not having the "rare" item/patterns -- the plan/pattern is generated considering long-term to globally reduce the bin usage.I suggest delete this line.}

In this paper, we propose to use the pricing in the dualism theory 
% \rjf{reference} 
to explicitly quantify the performance of different solution components, \eg, packing patterns for BPP, and use this information to guide the training of the algorithm to build up the solution. Intuitively, using a bin associated with a pattern can be seen as a way to satisfy a certain quantity of items. A group of bad patterns can lead to a shortage of bin supply for certain items. Therefore, patterns with low costs (waste) and high ability to satisfy overall demands will have higher prices.

% \rjf{Intuitively, these rare patterns will have a high price and the algorithm should pack the bins with the lowest packing price.} \rjf{check correctness. Here, you need to provide some simple intuitions how pricing works here.}  \myresponse{Contradictory on pricing concept. And, "rare" patterns is actually not the case.}
%The way that we quantify the quality of solution building blocks (\eg, packing patterns) is through the pricing in the dualism theory. 
The building blocks or patterns are dynamically generated by taking into account both the objective in Eq.~\eqref{eqn:packing} and constraints in Eq.~\eqref{eqn:demand} and Eq.~\eqref{eqn:constraint} through pricing. We describe this more in detail in the next section.

\subsection{Pricing and duality}
Duality is the principle that an optimisation problem could be viewed as two related problems with same data: the primal problem and dual problem. 
Consider the following standard formulation for optimisation problem (denoted as \textit{primal problem}) defined in Eq.~\eqref{eqn:primal_obj}-\eqref{eqn:primal_equal_const}~\citep{bonnans_convex_2019}, where $f$ is objective function, $u$ and $v$ are constraints, 
\begin{eqnarray}
\textit{Primal:~~~~~~} 
\mbox{minimise~~~} \quad    & f(\mathbf{x}) &  \label{eqn:primal_obj}   \\
\mbox{subject to:} \quad &  u_i(\mathbf{x}) \le 0 &  i=1,\cdots, m \label{eqn:primal_inequation_const}\\
& v_j(\mathbf{x})  = 0   &  j=1,\cdots, n\label{eqn:primal_equal_const}
\end{eqnarray}

We can then write the associated Lagrange function 
$L(\mathbf{x}, \mathbf{\lambda}, \mathbf{\eta}) = f(\mathbf{x}) + \sum_{i=1}^m\lambda_i u_{i}(\mathbf{x}) + \sum_{j=1}^n \eta_j v_j(\mathbf{x})$, 
where  $\lambda_i$ and $\eta_j$ are Lagrange multipliers. 
It is clear that finding an $\mathbf{x^*}$ that minimises $L(\mathbf{x}, \mathbf{\lambda}, \mathbf{\eta})$ with proper $(\mathbf{\lambda}, \mathbf{\eta})$ set can also get the optimal solution of primal problem. 
We denote $g(\mathbf{\lambda}, \mathbf{\eta}) = \inf_{\mathbf{x}} L(\mathbf{x}, \mathbf{\lambda}, \mathbf{\eta})$
as the \textit{dual problem} that aims to find a lower bound of the primal problem, where $\inf_{\mathbf{x}} L$ denotes the lower bound of $L$.  

% \begin{eqnarray}
% \textit{Dual:~~~~~~} 
% \mbox{maximise~~~} \quad    & g(\mathbf{\lambda}, \mathbf{\eta}) &  \label{eqn:dual_obj}   \\
% \mbox{subject to:} \quad &  \lambda_i \ge 0 &  i=1,\cdots, m \label{eqn:dual_inequation_const}
% \end{eqnarray}
% %
% & v_j(\mathbf{x})  = 0   &  j=1,\cdots, n\label{eqn:dual_equal_const}

Maximising the dual problem will obtain a set of Lagrange multiplier $(\mathbf{\lambda}^*,\eta^*)$ that identifies the effects that a certain constraint will have on the objective value. Such a value is also called \textit{shadow price} in economy and management community~\citep{KUOSMANEN2021666}. We extend the term \textit{price} to describe the process of evaluating key components in a candidate solution. Although the dual problem is also intractable computationally for most COPs, it becomes solvable when $f, u, v$ are linear functions, which is the case for the relaxed versions of many packing problems.  

The concept of duality has attracted much attention in many learning-related communities, with applications in dialogue~\citep{Jinpeng-Stylized-2021}, translation~\citep{wang2019MultiDualLearn}, etc. In offline combinatorial optimisation, there is rich literature that jointly applies dualism and pricing to solve large-scale mixed integer programming problems in a branch-and-price framework~\citep{DELORME20161}, which is essentially an iterative procedure 
% to repeatedly solve a restricted master problem (RMP) 
to repeatedly solve a dual problem
and a pricing sub-problem either exactly~\citep{wang_subgradientsimplexbasedcutting_2010} or approximately~\citep{xue_hybridpricingcutting_2021}. 

Existing research for applying duality for online problems is limited to heuristic analysis~\citep{gupta2020interior}. There lacks systematic research in pattern discovery and analyses in data mining community by leveraging the mathematical rigour of pattern generation and optimisation in operations research. In this paper, the proposed \texttt{CGPP} framework combines duality-based integer programming method and online pattern learning to generate high-quality patterns for packing incoming items. In addition to the benefits associated with patterns in terms of interpretability, the proposed approach could achieve superior solutions in terms of bin usage compared to other online algorithms and could, in some cases, achieve comparable performance to those offline approaches where all the items are known to the packing algorithm. 

\section{Literature review}
%\subsection{Algorithms for Bin Packing Problem}
\label{sec:literature_algo}

Bin packing problem has close connections to many real-world applications, \eg, memory management in modern computer architecture~\citep{Junjie-learning-2021}, healthcare management~\citep{abdalkareem_healthcareschedulingoptimization_2021c} and logistics~\citep{ALI2022108122}. It is probably one of the most studied combinatorial optimisation problems. Many research works have focused on approximate algorithms with provably guaranteed gaps to the optimal solution~\citep{gupta2020interior}. The most common ones are rule-based algorithms~\citep{coffman_binpackingapproximation_2013}, which deal with both online and offline BPP problems. \citet{coffman_binpackingapproximation_2013} provided a comprehensive review of classical bin packing heuristics. 

Recent solution methods for offline BPPs exploit the structural properties of BPP's integer programming formulations via exact algorithms like branch-and-bound schema~\citep{Zhang2020ColumnGA} and branch-and-price methods~\citep{DELLAMICO2020104825}, but the computational time varies significantly between instances and the methods are therefore not suitable for real-time decision making. 
% Metaheuristics/hyperheuristics/RL as black box optimisation
Another strand of research efforts is the data-driven based methods that exploit the distributional information of the random variables from the training data, including the use of genetic programming based hyper-heuristic to train a packing strategy/policy~\citep{5586388}, and the evolutionary algorithms for evolving rules to select the most appropriate packing heuristics at each decision point~\citep{lopez-camacho_unifiedhyperheuristicframework_2014}. 

Although the concept of patterns has been applied in offline combinatorial optimisation problems~\citep{LIU2021107175, Junjie-learning-2021,Dong-deep-2021}, it is not actively studied for online combinatorial optimisation problems until recently.
\citet{angelopoulos2022online} introduced \texttt{ProfilePack}, 
which utilises offline optimal solutions of a section of item sequence to generate a future packing plan that is used to guide the packing in real-time. Although the concept of pattern is not explicitly discussed in the paper, the high frequency packing examples in the offline optimal solution serve as templates to guide packing.
\citet{lin2024PatteAlgorFuzzy} developed another pattern-based packing method \texttt{PatternPack} for large-scale bin packing problem. It generates the pattern set by splitting the bin capacity into several fragments regardless of distribution of items while the plan is generated adaptively through statistical learning, coupled with a fuzzy logic enhanced pattern generation and selection strategy.
% which treats statistics from historical observation of item sequence as the prediction for incoming items, and generates the associated packing plan by heuristic and fuzzy packing strategy. 
These works demonstrate the potential of pattern-based methods for encoding and analysing the historical observation during packing. % with confirmed lower bound in terms of competitive ratio. 
However, these algorithms heavily rely on assumptions of simple stationery distributions, which may limit their performance for online COPs with non-stationary distributions.

Deep reinforcement learning (DRL) has gained growing attention in combinatorial optimisation~\citep{BENGIO2021405}, including 
%offline BPP~\citep{Cai2019ReinforcementLD,laterre_rankedrewardenabling_2018}, 
routing problems~\citep{Lu2020A} and graph-based problems~\citep{dai_learningcombinatorialoptimization_2018}. %\rjf{references, and list a few others with references.}
In most cases, BPP is formulated as a Markov Decision Process (MDP) through which uncertainties can be effectively handled by training on a large set of problem instances, and reinforcement learning methods are designed to tackle the problem in an end-to-end manner~\cite{zhang_attend2packbinpacking_2021}. For 1D bin packing problem, \citet{hubbs_orgymreinforcementlearning_2020} and \citet{ balaji_orlreinforcementlearning_2019} established a set of environments for classical operations research and associated DRL benchmarks. 
The benchmark end-to-end model  for 1D online bin packing problem proposed by \citet{balaji_orlreinforcementlearning_2019} is a simple multi-layer perceptron trained by PPO, achieved online packing by minimising the total Sum-of-Square potential~\citep{coffman_binpackingapproximation_2013}.
% aiming to minimising the total Sum-of-Square potential~\citep{coffman_binpackingapproximation_2013}. 
% Other researchers focused more on offline BPP. 
% \citet{Cai2019ReinforcementLD} utilised DRL to generate the initial solution for simulation annealing in order to solve offline 1D BPP. 
Additionally,
% \citet{laterre_rankedrewardenabling_2018} modelled the improvement of BPP solutions as a self-playing game, and applied Monte Carlo Tree Search to search a better solution. 
\citet{Junjie-learning-2021} developed SchedRL, a deep Q-learning method  with a specific reward design for online virtual machine scheduling, which can be modelled as an online variable-sized BPP. \citet{zhao_online3dbin_2021,zhao2022learning} investigated the online 3D bin packing with a robot arm for logistics, where the state is visually perceived through a deep neural network in~\citep{zhao_online3dbin_2021} and a graph neural network is designed to extract the position embedding of items in~\citep{zhao2022learning}. One noticeable drawback of these data-driven methods is their weak generalisation across unseen uncertainty distributions or non-stationary distributions. In this work, we establish a crucial link between the data mining and operations research and propose a novel pattern-based learning and optimisation method. The increased generalisation and performance enhancement of the proposed method is achieved through dynamic generation and resue of high-value patterns by explicitly exploiting the information of the uncertainty distributions. 
% DRL based heuristic search for bin packing problems are investigated recently.
% \citet{zhang2022DeepReinfLearn} solve 2D strip packing problem by DRL hyper-heuristics. \citet{romer2024MatheDiscoProgr} adapt a large language model to search and evolve heuristics for online bin packing problems, which is pretrained and fine-turned by RL.

% These works often achieve better performance comparing with heuristics,  
% \rjf{reference}
% and reasonable performance with lower time cost comparing with mixed integer programming solvers. 
% \rjf{reference}.\myresponse{It is a summary of previous works, do you mean add ref for methods they compared to?}  
%\rjf{the description is too brief.} %and \citet{balaji_orlreinforcementlearning_2019} 

% DRL agents have also been developed for online BPP. 

\section{Proposed pattern based method}
For online bin-packing problems, the constraints have as much impact on solving the optimisation problem as the optimisation objective which is often overlooked by most existing methods. 
We propose a general framework Column Generation Plan-and-Pack (\texttt{CGPP}) that adopts explicitly the dualism of COP to assist pattern based solution building. 
To do this, we first reformulate the BBP problem based the concept of patterns, then describe the key steps and modules in the CGPP framework, including mechanisms to handle uncertainty and imperfect distribution predictions. We explain how the dynamic pattern discovery through pricing could handle both the objective and the constraints well, leading to significant performance improvement for online BPP.

\subsection{Pattern based reformulation}
% \rjf{Motivations/existing challenges}

% \rjf{Process flow/ introduction of each building block.}

% \rjf{Discussion: Differences from existing solutions and benefits of proposed method.}

Formally, in our online BPP, we assume a problem instance as a finite sequence of items of length $N$, with index $i = 1,2,...,N$. Each item belongs to a finite type $t = 1,2,...,T$, which is associated with a size $s_t$. The quantity of item type $t$ in the sequence is defined as $q_t$, and its value is determined by sampling from a given distribution $D$.  %and its quanity, the ratio $\frac{q_t}{n} \sim D$.
In practice, the stochastic process of items could be more complex in the sense that the distribution could change over time. In this case, it becomes a non-stationary distribution problem which is harder to solve. Both stationary and non-stationary distributions of problems are studied in this paper. 
% We also investigate another type of problem in this research: the distribution is changing periodically.

We define a \textit{pattern} as a vector of the quantity of all item types that can be packed into a bin, $\mathbf{p}^h = (p_1^h, p_2^h,...,p_t^h,...,p_T^h)$ and $p_t^h=0$ means the item type $t$ does not appear in this pattern and $h$ is pattern index. %Therefore, a bin can be defined as a three-tuple: $b = (\mathbf{c}, \mathbf{p}, B)$, where $B$ is the capacity of bin. 
We denote $\mathbb{P}$ be the set of all feasible patterns. 
% Any $\mathcal{P}$ is a subset of all possible pattern set $\mathbb{P}$:
\begin{equation}
\mathbb{P}=\{\mathbf p^h | \sum_{t=1}^T{p_t^h s_t \le B}, p_t^h \in \mathbb{N}\} 
\end{equation}
% The total available pattern set is defined as $\mathcal{P} = \{\mathbf{p}_1, \mathbf{p}_2,...\mathbf{p}_m\}$ with associated pattern usage count $z_j, j =1,..,m$.
The original BPP formulation defined in Eq.~\eqref{eqn:objective}-\eqref{eqn:demand} can be re-formulated as follows:
\begin{eqnarray}
 \mbox{minimise} & ~\sum_{\mathbf{p}^h\in \mathbb{P}} z^h, &\label{eqn:pattern_bpp}\\
\mbox{subject to}& ~\sum_{\mathbf{p}^h\in \mathbb{P}} {p_t^h z^h  \ge q_t} &\quad \forall t=1,..,T\label{eqn:pattern_bpp_constr}
\end{eqnarray}
where $z^h\in \mathbb{N}$ is the decision variable, denoting the quantity of pattern $\mathbf{p}^h$ being used in a solution. 

In most cases, the feasible pattern set $\mathbb{P}$ is not prohibitively large and model Eq.~\eqref{eqn:pattern_bpp}-\eqref{eqn:pattern_bpp_constr} cannot be solved directly. In our method, the optimisation starts from a restricted pattern set $P \subset \mathbb{P}$ with total $m$ patterns. Then new patterns with potential to improve the objective value are iteratively added to the set $P$ by solving a pricing sub-problem until no solution-improving patterns can be found. 

\subsection{Framework overview}

\begin{figure*}[!htp]
\includegraphics[width=1\linewidth]{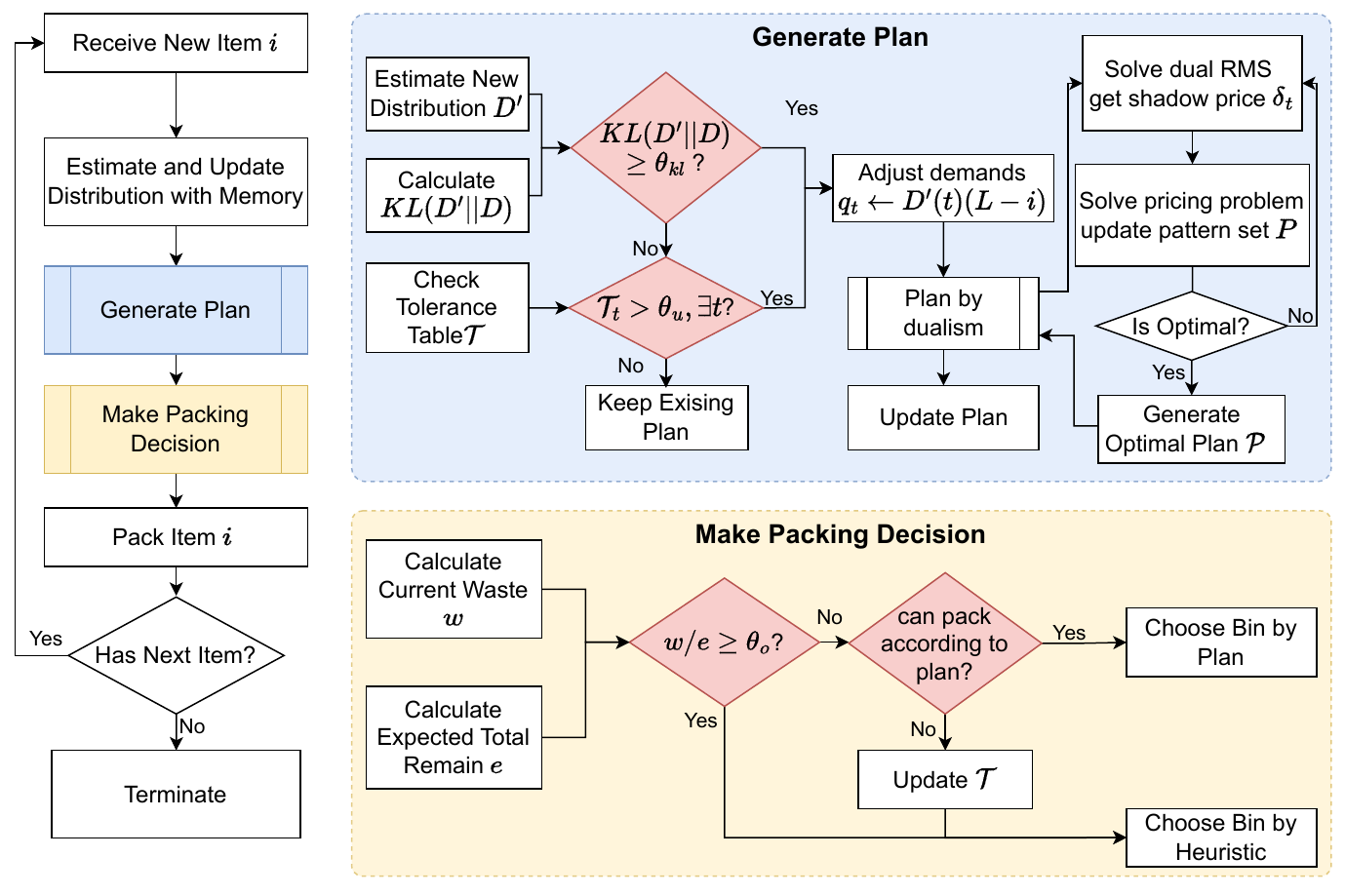}
\caption{The general framework of \texttt{CGPP}. Left: main procedure loop of iteratively packing items. Right: two critical module of the algorithm. Red rhombuses: uncertainty handling mechanism.}
\label{fig:framework}
\end{figure*}
The proposed \texttt{CGPP} framework is shown in Figure~\ref{fig:framework}.
The algorithm consists of three main stages: distribution estimation, plan generation and packing. The distribution estimation module aims to use a short-term sequence memory to estimate the real-time distribution of the random variables if the information is not known. The blue rectangle represents the plan generation procedure that applies dualism pricing method to identify good patterns and generates a plan to guide future packing. The yellow rectangle describes the adaptive packing process in CGPP with the guidance of the packing plan and the fallback strategies in the event of bad estimation errors. The details shall be described in the next few subsections.  %subsections~\ref{sec:planning}-\ref{sec:packing}.

\subsection{Planning}
\label{sec:planning}

As stated previously, obtaining full $\mathbb{P}$ is often not possible in most cases. Instead, the optimisation starts from a restricted master problem (RMP) formulated on a  subset $P \subset \mathbb{P}$ which guarantees a feasible solution but not the quality. A trivial way for the initial $P$ is to define a set of patterns in which each pattern packs one item type only. In the subsequent steps, the algorithm repeatedly generates new high-valued patterns to be added to $P$ and solves updated RMP, until no new pattern can be found to improve the solution further. The resulting solution by pattern frequency $z^h$ defines a packing \textit{plan}, $\mathcal{P}$. 
%A packing \textbf{plan} is then  formalised as $\mathcal P = (P, \mathbf{z}, \mathbf{e})$, where $\mathbf{z} = (z_1,..,z_m)$ represents the planned account of patterns and vector $\mathbf{e} = (e_1,...e_m)$ is used to track the execution status of patterns. 
By restricting to a small set of the patterns, we deal with a much smaller RMP and only high-valued patterns are added to the problem, considerably reducing the computational time. 

In order to identify high-valued patterns, we first obtain the shadow price $\delta_t$ for each constraint in Eq.\eqref{eqn:pattern_bpp_constr} through the dualism property. Then, the following sub-problem (called pricing problem or pattern generation) is solved:
\begin{eqnarray}
\mbox{minimise} &~ 1-\sum_{t=1}^T\delta_t p^*_t \label{eqn:pricing_obj}\\ 
\mbox{subject to } & ~ \sum_{t=1}^T s_t p^*_t \le B\label{eqn:pricing_constr}
\end{eqnarray}
and the resulting solution $\mathbf{p}^*=(p_1^*,p_2^*,...,p_T^*)$ defines a new pattern to be added to $P$. Eq. \eqref{eqn:pricing_constr} is the packing constraint. The pattern generation process stops once the objective value of Eq. (\ref{eqn:pricing_obj}) becomes non-negative which indicates that all potential cost-reducing patterns have successfully been discovered and resulting solution of RMP becomes optimal. The problem \eqref{eqn:pricing_obj}-\eqref{eqn:pricing_constr} is a knapsack problem that can be solved efficiently when the number of item types is not large which is the case for most real-world applications. 
%After the pattern identification, the primal problem is solved with the updated pattern set, through which, we obtain the optimal plan $\mathcal{P}$ and its associated quota $z_\mathbf{p}$ for each pattern $\mathbf{p}\in P$. $z_\mathbf{p}=0$ means the given pattern is not used in the plan $\mathcal{P}$.
In Algorithm~\ref{algo:planning}, lines 4-14 describe the pattern generation process.

Apart from the pattern set, another important component in RMP is forecasting the demand for each item type. We do this by dividing the whole item sequence into non-overlapping, equal-length sub-sequences (or sections), each of which is used to estimate the distributions of item types. Denote the section length to be $L$, we obtain a memory window of size $k \leq L$. At packing step $i$, a distribution $D'$ is estimated with observation of items from $i-k$ to $i$. We utilise Kernel Density Estimation (KDE) to determine the proportion of item types. This technique learns an appropriate linear combination of several Gaussian distributions, all sharing a same standard deviation but differing in their means. The model parameters are trained incrementally during the packing process, allowing the distribution estimation to adapt and improve over time. The demand of item type $t$, $q_t$, is set to the expected quantity of the type left in the remainder of the item sequence, i.e. $q_t = D(t)*(L-i)$.

% \change{TODO}

% \change{By solving relaxed RMS we obtain both the relaxed RMS objective and the shadow prices for each item type.
% When the relaxed RMS is solved, a key step is to iteratively find a good pattern to add into the pattern set $P$ to avoid enumerating all feasible patterns. To achieve this, a sub-problem (called pricing problem or pattern generation) is repeatedly solved to find a new promising pattern $\mathbf{p}*$ that minimises the reduced cost function (Eq. \eqref{eqn:pricing_obj}) for given shadow prices $\delta_t$ associated with demand constraint \eqref{eqn:pattern_bpp_constr}, while satisfying the packing constraint \eqref{eqn:pricing_constr} at the same time. }

% Algorithm~\ref{algo:planning} line 4-14 describes the general plan generation strategy.
%

\begin{algorithm}[!thbp]
\caption{Planning through Pricing by Dualism at item $i$}
\label{algo:planning}
\textbf{Input}: Memory length $k$, Previous plan $\mathcal{P}$, Priori distribution $D$, Underestimation tolerance table $\{\mathcal{T}_t\}, t=1,..,T$, Section length $L$\\
\textbf{Parameters}: Distribution threshold $\theta_{kl}$, Underestimate tolerance $\theta_{u}$ %, Max iteration $j_{max}$, 

\begin{algorithmic}[1]
    \STATE Estimate the current distribution $D'$ with item $i-k, .., i$
    \IF {$KL(D'||D) \geq \theta_{kl}$ \OR $ \exists \mathcal{T}_t \geq \theta_{u}, t =1,..,T$}
        \STATE {$D \leftarrow D'$}
        \STATE{Estimate remain demands $q_t \leftarrow D(t)(L-i), t=1,2,..,T$}
        % \STATE {Remove non-filled item slots in opened bins from $\mathbf{d}$} 
        % \STATE{Add sufficient count for each item $i$ to $\mathbf{d}$}
        \STATE{Initialise pattern set $P$}
        % \FOR{$j = 1,... j_{max}$}
        \WHILE{$1-\sum_{t=1}^T\delta_t p^*_t < 0$}
            \STATE{Solve Dual problem of RMP, obtain shadow prices $\delta_t$}
            \STATE{Solve model \eqref{eqn:pricing_obj}-\eqref{eqn:pricing_constr} with $\delta_t$, get $\mathbf p^*$}
            \STATE{Update pattern set $P \leftarrow P \cup \{\mathbf p^*\}$}
            % \IF {$1-\sum_{t=1}^T\delta_t p^*_t \leq 0$}
                % \STATE{Break}
            % \ENDIF
        % \ENDFOR
         \ENDWHILE
        \STATE{Solve integer programming model \eqref{eqn:pattern_bpp}-\eqref{eqn:pattern_bpp_constr} to get updated plan $\mathcal{P}'$}

        \STATE{$\mathcal{P} \leftarrow \mathcal{P}'$}
        \STATE{Clear $\mathcal{T}$}
    \ENDIF
    \RETURN{$\mathcal{P}$}

\end{algorithmic}
\end{algorithm}

\subsection{Plan-based packing} 
\label{sec:packing}

In an ideal world, the plan generated by Algorithm!\ref{algo:planning} is implemented exactly. However, because of forecast errors in demands, additional work is required during the actual packing (see Algorithm~\ref{algo:packing}). For a given packing plan $\mathcal{P}$, each newly opened bin is assigned to a pattern from the plan, implicitly specifying the type and quantity of items that should be packed into. Only items that match the assigned pattern can be packed in the corresponding bin. 
%The bin is therefore defined as a tuple $b=(\mathbf{r}, \mathbf{p})$, where $\mathbf{r}$ represents the actual item packing count vector and $\mathbf{p}$ represents the assigned pattern for the bin. We use the term ``match'' to describe whether a bin can pack an item according to its assigned plan. We say a bin matches an item of type $t$ if and only if $0 \leq r_t < p_t$. 
% When there is no matching bin,  an arbitrary pattern satisfying $p_t > 0$ and remain quota.
% Given input item with type $t$, a matched pattern means it has non-zero associated quantity $p_t$. 
% Obviously, when executing the plan, the used count of a matched pattern should not exceed it's planned quota. 
% Whether a bin $b$ matches an item is a little bit complicated: the bin should both have enough unfilled space and the quantity of $t$ already packed should not exceed the pattern allocated.

Upon arrival of an item of type $t$, the algorithm firstly packs it into a matched open bin via procedure $\mathtt{pack\_item}$. If no opened bin matches the considered item, a new bin is opened and an arbitrarily feasible pattern in the current plan $\mathcal{P}$ is assigned to it. The considered item is then packed to this new opened bin. This is done by procedure $\mathtt{open\_bin\_with\_pattern}$. Once a pattern in the plan has been assigned to a bin, its use frequency in the plan must be updated accordingly. Obviously, when executing the plan, the used count of a matched pattern should not exceed it's planned quota $z^h$ in the plan. %When the plan is perfect, the strategy can precisely execute the plan, resulting in a high-quality online solution. 

% In case the plan exhausted matched pattern or for uncertainty control purpose, the item will be packed with a fallback heuristic, e.g. best-fit in this research. Such fallback heuristic is executed by procedure $\mathtt{fallback\_pack}$ in the algorithm. Note the fallback decision can both assign an item into an existing bin regardless its pattern, and open a new bin with default pattern associated. In the former case, the incoming item will replace items in pattern but not packed with minimal size difference. In the latter case, a default pattern is assigned that contains only $t$ and $B-t$.
% ?The bins with empty pattern will only be filled by fallback heuristics.
When an item's demand is underestimated, at some point, there would be no feasible pattern available in the plan to pack this item. In such a case, the item is packed by a fallback heuristic, e.g. best-fit in this research. The fallback heuristic is executed by procedure $\mathtt{fallback\_pack}$ in the algorithm. Note that the fallback heuristic either assigns the item into an existing bin which would inevitably break its pattern requirements, or opens a new bin to pack the item. %In the former case, the incoming item will replace items in the pattern but not packed with minimal size difference. In the later case, a default pattern is assigned that contains only $t$ and $B-t$.
\begin{algorithm}[!t]
\caption{Pattern Based Packing Strategy at item $i$}
\label{algo:packing}
\textbf{Input}: Packing plan $\mathcal{P}$, Item distribution $D$, 
                Opened bins $\mathcal{B}$, Section length $L$ \\
\textbf{Parameters}: Overestimate tolerance threshold $\theta_{o}$ 

\begin{algorithmic}[1]
    \STATE{Calculate remain size $e \leftarrow (L - i) \sum_{t=1}^T s_t D(t)$}
    \STATE{Calculate total empty space of opened bins $w$}
    \IF{$w / e \geq \theta_{o}$}
        \STATE{$\mathtt{fallback\_pack}(i)$}
    \ELSE 
        \IF{ There exists an open bin $b$ whose pattern matches $i$}
            \STATE{$\mathtt{pack\_item}(i, b)$}
        \ELSIF{Pattern $\mathbf{p}$ in the plan matches $i$}
            \STATE{$b \leftarrow \mathtt{open\_bin\_with\_pattern}(\mathbf{p})$}
            \STATE{$\mathtt{pack\_item}(i, b)$}
        \ELSE
            \STATE{Update tolerance table $\mathcal{T}_{t}$ for item $i$}
             \STATE{$\mathtt{fallback\_pack}(i)$}
        \ENDIF
    \ENDIF

\end{algorithmic}
\end{algorithm}

\subsection{Uncertainty handling}
\label{sec:uncertainty}

Due to the stochastic nature of the problem and the imperfect estimation of quantities of items, a gap will always exist between the actual realisation of the problem instance and the forecast demands, leading to sub-optimal solutions. The challenge becomes greater when the distribution of item types is unknown and is subject to changes over time. 
%Two major factors contribute to the occurrence of uncertainty: the incomplete information nature of the problem, and the disparity between the estimated and actual item distributions. 

% The uncertainty caused by insufficient information can not be avoid in online problem for the exact item quantity is not available even for a section by definition.
% However, the bias between estimated and actual distributions might caused by systematic factors and can be fixed. For example, the item distribution changes dynamically, which is a common behaviour under real-world application\cite{9209730}. This will lead to wrong guess on item quantity demands, results in mismatch plan being generated. 
% \texttt{CGPP} checks whether such systematic error happens periodically at the beginning of each section. The distribution $D$ of last section is maintained, and starting $k$ items at current section are used to estimate new distribution $D'$. The Kullback-Leibler (KL) divergence $KL(D' || D)$ is utilised to compare the similarity between two distributions and a threshold $\theta_{kl}$ is applied. Violating the threshold means the past plan is not suitable for incoming items. Therefore, the planning procedure is executed again based on the newly estimated distribution to match the current observation.
The uncertainty caused by insufficient information cannot be avoided in online problems.
However, the gap between estimated and actual distributions might be caused by systematic factors which can be reduced. For example, when items' distribution changes dynamically over time, which is common in real-world applications \citep{9209730}, algorithms trained on stationary distributions could perform poorly. 

In our \texttt{CGPP} method, this challenge is dealt with by checking the distribution periodically at each section of the item sequence. The currently adopted distribution $D$ and the real-time distribution $D'$ estimated from the past $k$ items are compared using the Kullback-Leibler (KL) divergence $KL(D' || D)$. If the difference exceeds a predefined threshold $\theta_{kl}$, a new plan is generated by calling the planning procedure again based on the latest estimation of the distribution.

% The uncertainty will cause the estimated demand to be either underestimated or overestimated, hence introduce plan error. Underestimation happens when the actual quantity of items of a particular type is lower than estimated while overestimation arises when the actual number of items exceeds the estimation.
% If the algorithm do not handle uncertainty, the overestimated items tend to cause waste in space while the underestimated items will be packed through fallback strategy inefficiently. This results in numerous open bins waiting for items that will never arrive, or disrupts the predefined packing plan. 

The uncertainty can lead to errors in the estimated items' quantity, in the forms of either underestimation or overestimation.
Underestimation arises when the actual number of items exceeds the estimation, while overestimation occurs when the actual quantity of items of a specific type is lower than estimated.
Without special attention, the overestimated items are likely to result in wasted space, while the underestimated items will be packed inefficiently using fallback strategies. This can lead to the presence of numerous open bins waiting for items that will never arrive or disrupt the predefined packing plan.

% Underestimation indicates that the model might be overly pessimistic towards certain types. We maintain an uncertainty table $\mathcal{T}$, which keeps a record of the number of items not included in the current plan. The tolerance level for each item type is controlled by a threshold quantity $\theta_{u}$. This allows each item to be excluded from the plan for a maximum of $\theta_{u}$ times. Any $\mathcal{T}_t$ greater than $\theta_u$ indicates the expected demand of that type is far less than reality, hence the plan should be altered. We simply estimate the item distribution again as well as regenerate packing plan. Otherwise, an underestimated item is tolerated and packed by fallback heuristic.
%Underestimation indicates that the algorithm might be overly pessimistic towards demands of certain types. 
To address underestimation, we maintain an uncertainty table $\mathcal{\{T_t\}}$ that tracks the number of items not included in the current plan. The tolerance level for each item type is determined by a threshold quantity $\theta_{u}$. This means that each item can be excluded from the plan for a maximum of $\theta_{u}$ occurrences. If the count $\mathcal{T}_t$ for item type $t$ exceeds $\theta_u$, it indicates that the expected demand for that type is significantly lower than the actual demand, and the plan needs to be adjusted. In such cases, we re-estimate the item distribution and regenerate the packing plan accordingly. Otherwise, an underestimated item is tolerated and packed by fallback heuristic.

On the other hand, detecting overestimation is challenging until the very end of the item section. An opened but unfilled bin may be filled later with the planned items, or it could wait for an item that never comes according to plan, resulting in substantial waste of bin capacity.
To measure the risk of overestimation occurring, we introduce a risk ratio $w/e$, where $e$ is the expected sum of the remaining items in the section, and $w$ is the current total empty space across all opened bins. A higher ratio indicates a higher level of risks for following the current plan.
A threshold $\theta_{o}$ is introduced to express our tolerance for risk of overestimation. 
The bins are allowed to wait for incoming items until the threshold is exceeded, at which point, it becomes too risky to follow the plan, and the fallback heuristic is triggered to pack all remaining items instead.
%We allow at most $e\theta_o$ space of opened bins to be left empty, where $e$ is the expected sum of item remained in the section. The total empty space is defined as $w$ in Algorithm~\ref{algo:packing}. 

\section{Experimental results}
\label{sec:experiments}
We test the proposed method for a whole range of online BPP datasets with different characteristics in order to establish comprehensive evaluations and understanding of the strength and weaknesses of our method under different uncertainty conditions. More specifically, four distinct problem types are tested and details are given later in Sections~\ref{sec:exp-scope}-\ref{sec:experiments-large-scale}. Most experiments were set up with 20 instances, each having 20,000 items. Without explicitly stated otherwise, the bin capacity is set to 100, and the item sizes are in the range $[1, 100)$. 

We compare our method against \texttt{BestFit}, which is one of the most commonly used online heuristics due to its robustness across different scenarios and low competitive ratio (1.7), and three other state of the art methods for online BPP, namely \texttt{ORL} \citep{balaji_orlreinforcementlearning_2019}, \texttt{ProfilePack} (or ProfP for brevity) \citep{angelopoulos2022online}, and \texttt{PatternPack} (or PatnP for brevity) \citep{sym14071301}. An additional comparison with \texttt{PatternPack}'s updated version \texttt{FPP} \citep{lin2024PatteAlgorFuzzy} is given in subsection~\ref{sec:experiments-large-scale}.

\subsection{Algorithm configuration}
\label{sec:model-setup}
% We used the gap between the classic L2 lower bound \citep{martello_lowerboundsreduction_1990} and objective values from different algorithms to evaluate their performance. It is proven to be lower than any offline or online algorithm and achieves less than $1\%$ error with optimal value \citep{martello_knapsackproblemsalgorithms_1990}. 
% We used the gap between the classic L2 lower bound \citep{martello_lowerboundsreduction_1990} and the objective values obtained from different algorithms to evaluate their performance. This lower bound has been proven to be lower than any offline or online algorithm and achieves less than $1\%$ error compared to the optimal value \citep{martello_knapsackproblemsalgorithms_1990}.
We utilised the discrepancy between the classic L2 lower bound \citep{martello_lowerboundsreduction_1990} and the objective values obtained from various algorithms to assess their performance. This lower bound has been demonstrated to be less than $1\%$ from the optimal value.

%Additionally, it is easy to be calculated for large-scale problem instances for its $O(n^2)$ time complexity. We report the gap of total used bin between any tested algorithm against L2 lower bound to measure performance.
% To establish a benchmark for comparison, we build an oracle offline solution to each instance.
% % where the solver possesses exact quantity of the items to be packed. 
% The oracle simply execute the planning stage with exact item quantity and assumes the packing procedure perfectly follows the plan.
% This allowed us to measure the solution quality of the proposed method by calculating the gap between the number of bins used by the proposed method and the oracle solutions. 

% If not specified, the \texttt{CGPP} in this work is configured as follows. The fallback strategy is set to be one-step best fit heuristic. The section length is set to be $L=1000$ with memory length $k=250$. For the threshold parameters, the KL-Divergence threshold $\theta_{kl} =0.1$, underestimate tolerance threshold $\theta_u=5$ and overestimate threshold $\theta_o = 0.8$. 
If not specified, the \texttt{CGPP} in this work was configured as follows. The fallback strategy was set to be the one-step best fit heuristic. The section length was set to $L=1000$ with a memory length of $k=250$ based on some initial trials. For the threshold parameters, the KL-Divergence threshold $\theta_{kl} = 0.1$, the underestimate tolerance threshold $\theta_u = 5$, and the overestimate threshold $\theta_o = 0.8$.
% The adaptive version of \texttt{ProfilePack} \cite{angel2023OnlinBinPacki} is used in this experiment. 
% The adaptive \texttt{ProfilePack} is based on the hybrid profile packing strategy, which mix the profile packing and first fit heuristic by some ratio. We set the hybrid ratio $\lambda_{PP}=0.5$ as it is been reported to have most stable performance against large prediction error\citep{angelopoulos2022online}. 
The hybrid \texttt{ProfilePack} algorithm was set up with parameter $\lambda_{PP}=0.5$ as suggested by \citet{angelopoulos2022online} since a low-error profile was not assumed in our experiment. 
On the other hand, \texttt{PatternPack} was configured with the same parameters reported in \citet{sym14071301}.
% \texttt{PatternPack} generates the pattern set by splitting the bin capacity into several fragments regardless of distribution while the plan is generated adaptively by statistical learning the distribution. \citet{sym14071301} defines the number of fragments $F_{PaP}$ and total number of patterns $P_{PaP}$ to control the pattern generation behaviour. 
% We configure the algorithm with $F_{PaP}=17$ and $P_{PaP}=100$ the same as \citet{sym14071301} reported. 
% Both \texttt{ProfilePack} and \texttt{PatternPack} applied a statistical learning approach to adaptively learn the problem distribution. This approach kept a sliding memory window and the frequency of item in the window is applied to calculate the estimated probability. 
% We set up the memory window length of two algorithms to be $k_{PP}=k_{PaP}=500$.
% Both \texttt{ProfilePack} and \texttt{PatternPack} applied a statistical learning approach to adaptively learn the problem distribution. This approach involved maintaining a sliding memory window, where the frequency of items in the window was used to calculate the estimated probability.
% We set up the memory window length for both algorithms to be $k_{PrP}=k_{PaP}=500$.
Both \texttt{ProfilePack} and \texttt{PatternPack} employed a statistical learning approach to dynamically learn the problem distribution. This approach entailed maintaining a sliding memory window, in which the item frequencies within the window were utilised to estimate the probabilities.
For both algorithms, the length of the memory window to be $k_{PrP}=k_{PaP}=500$, same as the settings reported in the papers.
% Additionally, in Section~\ref{sec:experiments-large-scale} we also compared with a improved method  \citep{lin2024PatteAlgorFuzzy}
% To maintain a similar learning ability, 

% Similarly, \texttt{ProfilePack} also maintains a sliding memory window to evolve the distribution and generate profile for future. We set the window length $k_{pp}=500$ same as \texttt{CGPP}'s.

% The \texttt{ORL} benchmark was trained under a uniform distribution given items and capacity same with problem definition, if not mentioned specifically. We utilised the standard PPO algorithm~\citep{DBLP:journals/corr/SchulmanWDRK17} implemented by Tianshou framework~\citep{tianshou} to train the model, employing 100 parallel training environments. We adapt a similar environment as the BPP environment by \citet{balaji_orlreinforcementlearning_2019} with different distribution and item type settings according to problem. Each model underwent 500 epochs of training, which took approximately 600 minutes to complete. The policy net is a 3-layer ANN with hidden layer of size 256, while other model parameters kept same with \citet{balaji_orlreinforcementlearning_2019}. 

The \texttt{ORL} method in \citet{balaji_orlreinforcementlearning_2019} was re-implemented with the same reported settings. The algorithm used the standard Proximal Policy Optimisation (PPO) algorithm, with a 3-layer policy network and a hidden layer of 256 nodes. The model was trained on a uniform distribution set, where the items and bin capacity were the same as the problem definition, unless otherwise specified. It underwent 500 epochs of training, which took approximately 600 minutes to complete on our machine. 

All experiments were performed on a PC with an Intel Xeon Gold 6248R Processor with 24 cores and 48 threads, along with an Nvidia GeForce RTX 3090 graphics card.

\subsection{Experiments on different items' distributions}
\label{sec:exp-scope}

% This experiment set was designed to identify what type of distribution that \texttt{CGPP} performs good. In total 8 experiments are set up with uniform or normal distribution. We set up three derived experiments with same distribution but different item range configurations. The suffix B refers to biased distribution, with item size in range$[10, 60)$, while suffix S refers to symmetric distribution, with item size in range $[25,75)$. Specifically, for Normal-B, the mean of distribution is set to be $\mu=35$. The experiment wit suffix C refers to coarse experiment, with item in range $\{10, 20, ..., 90\}$. 

In order to evaluate the performance differences across different distributions by all algorithms, a total 8 datasets were set up with uniform or normal distributions as bases. They are named \textit{Uniform}, \textit{Normal}, \textit{Uniform-B}, \textit{Uniform-C}, \textit{Uniform-D}, \textit{Normal-B}, \textit{Normal-S}, and \textit{Normal-C}, respectively. Among them, 6 are derived datasets with a same distribution but different item range configurations. The suffix B refers to biased distribution, with item size in range $[10, 60)$, while suffix S refers to symmetric distribution, with item size in range $[25,75)$. Specifically, for Normal-B, the mean of distribution was set to be $\mu=35$ with the range same as Uniform-B. The experiment with suffix C refers to coarse experiment, with item sizes from set $\{10, 20, ..., 90\}$. 

Table~\ref{tab:bpp-scope} provides the results of this experiment set. It can be seen that \texttt{CGPP} outperformed other methods in most experiments, except for two of the symmetrically distributed sets (Normal and Normal-S), for which \texttt{BestFit} outperformed all other methods. %in ordinary Binomial and Binomial-S, where the item distribution is symmetric to 50, i.e. half of the capacity. However, the \texttt{CGPP} still generated solution close to Best Fit and better than other three algorithms on those two datasets. 
%Results showed that \texttt{CGPP} is effective across most distributions, including problems with non-symmetric distribution, or with relatively fewer item types. 
The proposed CGPP method tends to perform particularly well for uniform distribution instances. The performance gain against the second best method for these instances ranges between 17\% to 62\%. 
% Intuitively, an item of size $i$ in the symmetric distribution sequence 

\begin{table}[htbp]
    \centering

    \begin{tabular}{|l|c|c|c|c|c|}
    \hline

    Distribution & BestFit & ORL & ProfP &PatnP &CGPP\\

    \hline
    1. Uniform & 343.60  & 700.35 & 379.4 & 343.10 &\textbf{271.75}\\
    2. Uniform-B & 199.8  & 339.6 &  546.35 &  249.05 &\textbf{76.60}\\
    3. Uniform-S & 100.55  & 252.55 & 323.45 & 159.85 &\textbf{79.60}\\
    4. Uniform-C & 48.15  & 1186.9 & 1010.5 & 562.35 & \textbf{39.85} \\

    5. Normal & \textbf{490.9}  & 1591.5 &  1722.45 &2319.4 &507.55\\
    6. Normal-B &1012.65  & 1202.95 & 1165.85 &1012.65 &\textbf{954.70}\\
    7. Normal-S & \textbf{77.3}  & 1280.55 & 1205.60 & 1738.60 & 167.2\\
    8. Normal-C & 490.95  & 494.1 & 2290.40 & 494.0 &\textbf{489.65}\\
    % 5. Normal & \textbf{121.75}  & 1015.3 &  710.0 &595.55 &324.55\\
    % 6. Normal-B &282.5  & 342.35 & 559.0 &319.85 &\textbf{258.65}\\
    % 7. Normal-S & \textbf{90.05}  & 1140.15 & 1352.45 & 1231.2 & 400.2\\
    % 8. Normal-C & 473.05  & 478.9 & 1644.25 & 476.75 &\textbf{472.15}\\
    \hline
    \textbf{Overall average} & 345.5  & 881.1 & 1080.5 & 859.9 &\textbf{323.4}\\
    \hline
    \end{tabular}
    \caption{The average objective gaps to L2 bound by different algorithms for problems with different uniform and normal distributions.  Bold text represents best average results. ProfP refers to \texttt{ProfilePack} and PatnP refers to \texttt{PatternPack}.}
    \label{tab:bpp-scope}
\end{table}

\subsection{Packing with prior knowledge}
\label{sec:experiment-learning-cost}
% This experiment set was built to investigate whether a good priori can contribute to finding a good solution. Among the algorithm we discussed, \texttt{BestFit} does not rely on any learning mechanism; \texttt{ProfilePack} and \texttt{PatternPack} both apply statistical approach to gain the distribution information without any priori knowledge. On the other hand, \texttt{ORL} can be viewed as implicitly encode the distribution information by choose the training environment. In this experiments, \texttt{ORL} was trained on the same distribution as experiment. We also report the result of \texttt{CGPP} with exact distribution given as priori as \texttt{CGPP-L}.
This experiment set was built to investigate whether a good prior knowledge on distribution can contribute to finding a good solution. Among the algorithms we discussed, \texttt{BestFit} does not rely on any learning mechanism, while both \texttt{ProfilePack} and \texttt{PatternPack} apply a statistical approach to gain distribution information without any prior knowledge. On the other hand, \texttt{ORL} can be viewed as implicitly encoding the distribution information by choosing the training and testing datasets. In the experiments in this section, \texttt{ORL} was trained on the same distribution as the test datasets' distribution. Additionally, we report the results of \texttt{CGPP} with the exact distribution given as prior knowledge, referred as \texttt{CGPP-L}.

% To build a convincing comparison basis with \texttt{ORL}, we adapt three distributions, namely BW1, LW1, and PP1, as proposed by \citet{balaji_orlreinforcementlearning_2019}. These distributions has expected waste of $\Theta(1)$, $\Theta(\sqrt{n})$ and $\Theta(n)$. The bin capacity in experiment 9-311 was set to 9 while in 12-14 was 100, also following the configuration by \citet{balaji_orlreinforcementlearning_2019}. The number of experiment instances and associated item number are kept same as Section~\ref{sec:exp-scope}.

To establish a convincing comparison with \texttt{ORL}, we adopted three distributions: BW1, LW1, and PP1, as proposed by \citet{balaji_orlreinforcementlearning_2019}. These distributions have expected waste of $\Theta(1)$, $\Theta(\sqrt{n})$, and $\Theta(n)$, respectively. A total of 6 datasets are created (see Table~\ref{tab:bpp-priori}). In the first 3 datasets (BW1-9, LW1-9, PP1-9) the bin capacity was set to 9, while in datasets BW1-100, LW1-100 and PP1-100, it was set to 100, following the configuration by \citet{balaji_orlreinforcementlearning_2019}. The number of experiment instances and the associated number of items remained the same as in Section~\ref{sec:exp-scope}.
\begin{table}[htbp]
    \centering

    \begin{tabular}{|l|c|c|c|c|c|c|}
    \hline

    Distribution & BestFit & ORL & ProfP &PatnP &CGPP & CGPP-L\\

    \hline
    9. BW1-9 & \textbf{0.00} & \textbf{0.00} & 527.55&\textbf{0.00} & 0.25 & \textbf{0.00}\\
    10. LW1-9 & 103.60 & 156.50 & 388.70 & 103.60 & 223.55 &\textbf{101.45}\\
    11. PP1-9 & 154.90 & 472.75& 1187.25 &154.90 & 146.30 & \textbf{145.40}\\
    12. BW1-100 & \textbf{14.10} & \textbf{14.10} & 428.10&\textbf{14.10} &\textbf{14.10} & \textbf{14.10}\\
    13. LW1-100 & \textbf{0.00} &  \textbf{0.00} &  792.95&\textbf{0.00}&\textbf{0.00} & \textbf{0.00}\\
    14. PP1-100 & \textbf{0.00} &\textbf{0.00} & 702.15&\textbf{0.00}&\textbf{0.00} & \textbf{0.00}\\
    \hline
    \end{tabular}
    \caption{Experiment results on distributions proposed by \citet{balaji_orlreinforcementlearning_2019}. Bold text represents results with best average bin gap. ProfP refers to \texttt{ProfilePack} and PatnP refers to \texttt{PatternPack}. }
    \label{tab:bpp-priori}
\end{table}

% Table~\ref{tab:bpp-priori} represents the experiment results. For experiment 1 and 4-6, almost all models were able to achieve optimal, except \texttt{ProfilePack}. For experiment 2 \texttt{CGPP-L} outperformed \texttt{BestFit} and \texttt{PatternPack}, but \texttt{CGPP} failed. For experiment 3, \texttt{CGPP-L} also outperformed \texttt{CGPP}. These observation shows the learning distribution online will cause extra cost comparing with a good distribution guess. \texttt{ORL} is more unstable: it performs badly on experiment 3 even if it was trained on that distribution. On the other hand, \texttt{PatternPack} achieved good online learning strategy as the result was almost same as \texttt{BestFit} while \texttt{ProfilePack} achieves worst performance.
Table~\ref{tab:bpp-priori} shows the experiment results. For datasets 9, 12, 13 and 14, almost all methods were able to achieve optimal solution, except \texttt{ProfilePack}. For datasets 10 and 11, \texttt{CGPP-L} outperformed all methods. However, \texttt{CGPP} failed to obtain competitive results for dataset 10, indicating a potential weakness of the proposed method and the importance of utilising prior knowledge if available. It appears that \texttt{CGPP}'s online distribution learning mechanism misled the packing because its prior-knowledge version performed best. \texttt{ORL} has mixed performance: it performs badly on dataset 11 even after it was trained on that same distribution. On the other hand, \texttt{PatternPack} achieved good online learning strategy for this group of datasets as the result was almost same as \texttt{BestFit} while \texttt{ProfilePack} achieves worst performance. Overall, this group of datasets seem to be rather friendly for best fit which does not have learning. 

\subsection{Experiments on more complex distributions}
\label{sec:exp-complex}

% This experiment set investigated the performance on complex distributions. Two types of experiments were set up: the first group adapt dual-normal distribution from \citet{burke_providingmemorymechanism_2010}, as typical mixed distribution. The dataset contains both single normal distribution and dual normal distribution sets. We focused on the dual part, namely Burke 4-11, as experiment 15-22 represented. Each experiment contains 20 instances with each instance having 5000 items. 

This section assesses the performance of our algorithm on more complex distributions. Two groups of datasets were used. The first group adopted a dual-normal distribution suggested by \citet{burke_providingmemorymechanism_2010}, which is a typical mixed distribution. The datasets include both single normal distribution and dual normal distributions. Our focus was on the dual part, specifically Burke 4-11, through experiments 15-22. Each experiment consists of 20 instances, with each instance containing 5000 items.

% The second group investigates the performance when the distribution changed periodically. All three experiments shared same item size range and bin capacity configurations. The whole item sequence was divided into several equal-sized sections, each section was sampled from an independent distribution. We utilised two binomial distribution groups based periodic experiments: the Binomial-PS samples from Binomial distribution with $p=\{0.2, 0.35, ..., 0.7\}$ and Binomial-PB samples with $p=\{0.2, 0.3, ..., 0.6\}$. We also add a Poisson distribution group with $\lambda$ changed in range $\{5, 15, ..., 45\}$. For each instance, the section step was set to be 2000, with total 10 sections. Each instance repeatedly sampled from the given distribution group. 
The second group of experiments investigated the performance when the distribution changes periodically. All three experiments shared the same item size range and bin capacity configurations. The entire item sequence was divided into several equal-sized sections, and each section was sampled from an independent distribution. We utilised two groups of binomial distributions for the periodic experiments: Binomial-PS, which samples from a binomial distribution with $p=\{0.2, 0.35, ..., 0.7\}$, and Binomial-PB, which samples from a binomial distribution with $p=\{0.2, 0.3, ..., 0.6\}$. Additionally, we included a Poisson distribution group with its parameter varying in the set $\{5, 15, ..., 45\}$.
For each instance, the section size was set to be 2000, resulting in a total of 10 sections. %Each instance repeatedly sampled from the given distribution group.
% To assess the model's ability to recognise unknown distributions, the model is initialised with uniform distribution. The 

% Our method significantly outperforms the other approaches.

\begin{table}[!t]
    \centering

    \begin{tabular}{|l|c|c|c|c|c|}
    \hline
        Distribution & BestFit & ORL & ProfP &PatnP & CGPP \\

    \hline
    15. Burke-4 & 205.00   & 219.60 & 2260.0 & 178.65&\textbf{115.75}  \\
    16. Burke-5 & 167.50  & 179.55 & 202.05 &165.25 &\textbf{82.35} \\
    17. Burke-6 & 75.20   & 85.20 & 196.85 &115.45  &\textbf{53.50}  \\
    18. Burke-7 & 50.60 & 56.90 & 120.50 &78.90  &\textbf{37.70} \\
    19. Burke-8 & 180.55 & 209.20 & 198.4&172.00  &\textbf{102.9}\\
    20. Burke-9 & 145.55 & 157.25 & 165.6 & 140.25  &\textbf{80.25}\\
    21. Burke-10 & 96.55 & 104.20 & 215.7 & 98.10  &\textbf{57.95}\\
    22. Burke-11 & 54.55 & 61.05 & 185.15 &73.25 & \textbf{42.75}\\
    \hline
    23. Binomial-PS & \textbf{1430.5}  & 1703.3 & 2349.1 & 1712.0 &1437.3 \\
    24. Binomial-PB & 1310.9  & 1319.2 & 1551.7 & 1437.7 &\textbf{1302.6} \\
    25. Poisson & 202.9  & 235.5 & 738.0 & 204.9 &\textbf{177.9} \\
    % 9. Uni+Nor & 66.00 & 150.15 & 45.05 &45.72 &\textbf{14.68} \\
    % 10. Uni+Poi & 15.59 & 41.26 & 25.28 &54.98 &\textbf{10.53}\\
    % 11. Nor+Poi & 43.3 & 99.85 &37.59 & 29.07 &\textbf{7.63}\\
    % 12. Poi+Wei & 13.39 & 47.54 &37.59 & 29.07 &\textbf{7.63}\\
    \hline
    \end{tabular}
    \caption{Experiment results for dual distributions and periodic distributions measured by average bin gap to the L2 bound. Bold text represents best objective values.}
    \label{tab:bpp-complex}

\end{table}

% Table~\ref{tab:bpp-complex} illustrates the results on two types of experiments. For dual distribution set, \texttt{CGPP} outperformed other methods significantly. For periodic distributions, \texttt{CGPP} outperform others on experiment 24, 25 but got close result with \texttt{BestFit} on experiment 23.
Table~\ref{tab:bpp-complex} shows the results for the two types of experiments. In the dual distribution set, \texttt{CGPP} outperformed the other methods significantly. Compared with \texttt{BestFit}, the reduction in the gap to L2 ranges from 21.5\% to 50.1\%. Compared with \texttt{PatternPack}, the reduction is between 35.2\% and 52.2\%. 

In the periodic distribution experiments, \texttt{CGPP} performed similarly to \texttt{BestFit} for two Binomial distributions but gained clearly advantage for Poisson distribution. Compared to \texttt{ProfilePack} and \texttt{PatternPack}, \texttt{CGPP} again has significant advantages. 
%the other methods in experiments 24 and 25, but achieved similar results to \texttt{BestFit} in experiment 23.

\subsection{Experiments on large-scale Weibull distribution}
\label{sec:experiments-large-scale}
% To cover the item shape properly, the scale parameter was set to be $sc=50$.
% Weibull distribution \citep{casti2012WeibuBenchBin} is a family of probabilistic distribution that connects closely with bin packing applications such as VM management. This section is aiming to investigate the effectiveness on the Weibull distribution family. We established five kinds of weibull distribution with shape parameter $sh=\{0.5, 1.0, 1.5, 2.0, 5.0\}$.  Each experiment contained 5 instances with $10^6$ items. Additionally, we also adopted a periodic Weibull distribution group with section size of 4000 and shape parameter iterated from the 5 distributions above. Several model parameters are modified: for \texttt{CGPP}, the memory length is $k=1000$ and section length was $k=4000$, the underestimate tolerance was set to be $\theta_u = 20$. The memory window of \texttt{ProfilePack} and \texttt{PatternPack} were set to be $k_{PP}=k_{PaP}=1000$.

In this section, we aim to investigate the effectiveness of the Weibull distribution family, which is closely connected to bin packing applications such as VM management \citep{casti2012WeibuBenchBin}. We established five different Weibull distributions with shape parameters $sh=\{0.5, 1.0, 1.5, 2.0, 5.0\}$. Each experiment consisted of 5 instances with $10^5$ items.

In addition, we generated a group of datasets with periodic Weibull distributions, with the shape parameters shifting to the next one stated in the list above for every 4000 items (section size).

Some of the parameters were modified for this experiment in order to adapt to instances with very large number of items. For \texttt{CGPP}, the memory length was set to $k=1000$, the section length was set to $k=4000$, and the underestimate tolerance is $\theta_u=20$. Overestimate tolerance is $\theta_o=1.5$ since the sequence is long enough to pack items according to plan. The memory window of \texttt{ProfilePack}, \texttt{PatternPack} and \texttt{FPP} were set to be $1000$.

\begin{table}[!t]
    \centering

    \begin{tabular}{|l|c|c|c|c|c|c|}
    \hline
        Distribution & BestFit & ORL & ProfP &PatnP &FPP& CGPP \\

    \hline
    26. $sh=0.5$ &  \textbf{0.2}  & 736.2 &  11695.8 & \textbf{0.2} & 154.2&476.2   \\
    27. $sh=1.0$ & 153.8  & 1541.4 &   2110.2 & 134.4 & 219.4 &\textbf{82.2}  \\
    28. $sh=1.5$ &  608.6  & 2386.0 &  1272.6 &  811.4 & 349.2  &\textbf{94.6}  \\
    29. $sh=2.0$ &  1039.8  & 3098.6 & 906.8 & 1316.8 &477.6 &\textbf{133.8} \\
    30. $sh=5.0$ &  2150.6  &  2981.2 & 1641.6 & 2515.4 & 892.4 & \textbf{384.2}\\%\textbf{736.2} \\
    31. Periodic &  465.4 &  802.2&  2306.0 & 609.2 &259.4& \textbf{208.4} \\
    \hline
    \end{tabular}
    \caption{Experiment results for large-scale Weibull distributions measured by average bin gap to the L2 bound. Bold text represents best results.}
    \label{tab:bpp-large-scale}

\end{table}

% Table~\ref{tab:bpp-large-scale} represents the experiment results. \texttt{CGPP} significantly outperformed other methods in most cases when $sh\geq 1.5$, and outperformed \texttt{BestFit} when $sh=1.0$. \texttt{CGPP} performs badly when $sh=0.5$. A possible guess is due to the small-size item was involved heavily in the sequence, but the algorithm kept assume large-size items would come.                
Table~\ref{tab:bpp-large-scale} shows the experiment results. Again, \texttt{CGPP} outperformed the other methods for almost all datasets except when $sh=0.5$. One possible explanation for this is that the sequence heavily involved small-size items, but the algorithm continued to assume that large-size items would come in the future. The proposed method relies on a good forecast of both the type of items and their distributions. When the uncertainty is extremely high, it is probably better to revert to more myopic methods like best fit. 

For this group of experiments, we also included the results from a very recent algorithm FPP. As was shown in Table~\ref{tab:bpp-large-scale}, compared to its previous version PatternPack, FPP obtains mixed results. It did quite well for instances with $sh=1.5, 2.0, 5.0$ and periodic instances, obtaining the second-best results among all compared algorithm. However, it is outperformed by PatternPack for the instances with $sh=0.5, 1.0$, suggesting some robustness issues. 
% As a result, \texttt{CGPP} struggled to adapt to the characteristics of the distribution and find the optimal bin size.
% This might due to the item in this distribution was quite small, but the algorithm tend to assume large size item will come in future. 

\section{Discussion and analysis}

\subsection{Solution quality analysis}

Although the performance of algorithms is primarily measured by bin usage, we use filled rate of all opened bins to further investigate the solution quality and the packing process by different methods. The bin filled rate is defined as the percentage of the total size of items in a bin to its capacity.

We use two typical solutions from Uniform-B and Normal-B datasets in Section~\ref{sec:exp-scope} for analysis. Additionally, we analyse the results on the periodic Weibull dataset (experiment 31) to observe how different methods behave when faced with changing distributions.
% as Figure~\ref{fig:solvis-uni} and Figure~-\ref{fig:solvis-norm} illustrate the averaged filled rate by bin series. 

\begin{figure}[t]
\centering
    % \subfloat[Uniform-B]{ \includegraphics[width=.8\textwidth]{uni.pdf}}\\
    % \subfloat[Normal-B]{ \includegraphics[width=.8\textwidth]{normal.pdf}}\\
    \includegraphics[width=\textwidth]{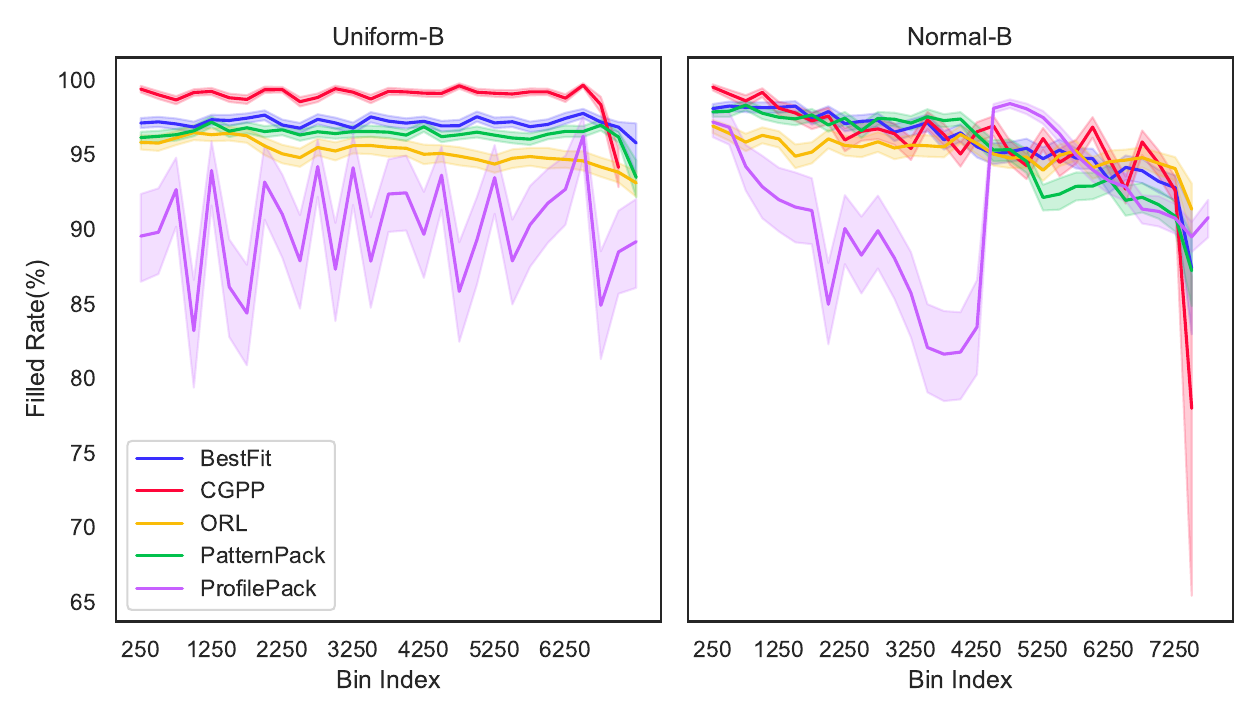}
    \caption{Average bin filled rate with confidence interval}
    \label{fig:solvis-combined}
    % \begin{subfigure}[b]{0.4\textwidth}
    %     \centering
    %     \includegraphics[width=.8\textwidth]{uni.pdf}
    %     \caption{Average Bin Filled Rate for Uniform-B}
    %     \label{fig:solvis-uni}
    % \end{subfigure}
    % \hfill
    % \begin{subfigure}[b]{0.4\textwidth}
    %     \centering
    %     \includegraphics[width=.8\textwidth]{normal.pdf}
    %     \caption{Average Bin Filled Rate for Normal-B}
    %     \label{fig:solvis-norm}
    % \end{subfigure}
\end{figure}

Figure~\ref{fig:solvis-combined} shows the filled rates of the bin index for the Uniform-B and Normal-B dataset. The bin series are arranged in the order of their opening steps. 
% All five bin solutions from the tested models are represented by polylines with different colours. The solutions generated by different online algorithms 
The polylines of different colours are used to illustrate the average fill rates of bins in the solutions generated by different algorithms for the given dataset.
% To provide a clearer representation, the average filled rate is calculated every 250 bins and represented by dark polylines. The 95\% confidence interval is also reported and shown as associated light areas.
The surrounded light areas of each polyline represent 95\% confidence interval. 

% On Uniform-B dataset, it is evident that \texttt{CGPP} achieves nearly full packing until the very end of the solution, significantly outperforming the other models. \texttt{BestFit}, \texttt{PatternPack}, and \texttt{ORL} perform similarly with filled rates ranging from 95\% to 97.5\%. T
It can be seen that both \texttt{PatternPack} and \texttt{ORL} is able to achieve filled rate over 95\% until at very late stages.
\texttt{BestFit} achieved slightly better filled rate than \texttt{PatternPack} and \texttt{ORL}, with average 97.5\%. 
The filled rates of \texttt{ProfilePack} are consistently worse than all other methods. The extreme fluctuations observed in \texttt{ProfilePack} also illustrates the algorithm's poor robustness. This phenomenon is likely due to \texttt{ProfilePack} lacking mechanisms for handling overestimation uncertainty, which is identified as the most wasteful resource, as discussed in Section~\ref{sec:uncertainty}.

% \begin{figure}[t]
% \centering
% \includegraphics[width=.8\textwidth]{normal.pdf}
% \caption{Average Bin Filled Rate for Normal-B}
% \label{fig:solvis-norm}
% \end{figure}
% Clearly the fill rate of CGPP is the highest at the beginning (albeit with fluctuations caused by not using the moving average), this could be caused by poor prediction. If you have perfect distribution knowledge, it would still outperform best fit. The fill rate of last few bins by CGPP also highlight the main drawback of CGPP method and could be an area for future improvement.
% Figure~\ref{fig:solvis-norm} represents the filled rates of the Normal-B dataset.
The filled rates of \texttt{BestFit}, \texttt{CGPP}, and \texttt{PatternPack} exhibited a decreasing trend on Normal-B dataset. Clearly the filled rate of \texttt{CGPP} is the highest at the beginning, albeit with fluctuations during the whole packing stage, which could be caused by imperfect prediction. The fast drop of filled rate of last few bins by CGPP also highlight one of the main drawback of CGPP method that the overestimation is not avoidable. Similarly, \texttt{PatternPack} also suffered with fluctuations and fast-drop by poor prediction.  %and could be an area for future improvement
The filled rate of \texttt{ProfilePack} experiences a significant drop in the middle of the bin sequence. This is likely due to poor predicted profile misguided the packing. \citet{angelopoulos2022online} claimed hybrid \texttt{ProfilePack} forced to pack items separately using online heuristic when the quantity of items packed following profile guidance reached to a threshold. This resulted in huge waste by not filling the space reserved for overestimated items.
% resulted huge waste by not filled overestimation.
% In the second half, \texttt{ProfilePack} forced to utilise the first-fit strategy \citep{angelopoulos2022online}, resulted huge waste by not filled overestimation.
% This is likely due to \texttt{ProfilePack} being forced to utilise the first-fit strategy in the second half of the bin sequence, which leads to wasting the overestimated items from the previous stage.
This further highlights the importance of addressing overestimation in the pattern based packing process. Therefore, eliminating the effect of poor prediction could be an area for future improvement for all prediction-based online algorithms. On the contrary, \texttt{ORL} behaved conservatively by maintaining most bins at similar level of filled rates for both datasets. Since \texttt{ORL} was trained on uniform distribution, such conservative strategy indicates it cannot generalise to other distributions.
% This shows \texttt{ORL} learned a conservative strategy, regardless of the distribution. As a consequence, it generated large amount of sub-optimal bins with guaranteed waste.
% However, this strategy proves to be sub-optimal as it encourages bins to hover around a 95\% filled rate, resulting in a guaranteed 5\% waste.
% This similarity may be attributed to the fact that both \texttt{CGPP} and \texttt{PatternPack} employ the best-fit strategy as a fallback strategy. The best-fit strategy aims to find the first bin with the least waste after packing, which allows bins opened earlier to be packed more effectively. Additionally, the limited quantity of small-sized items in Normal-B may not satisfy the bins that are waiting for them, resulting in potential waste.

\begin{figure}[t]
\centering
\includegraphics[width=.8\textwidth]{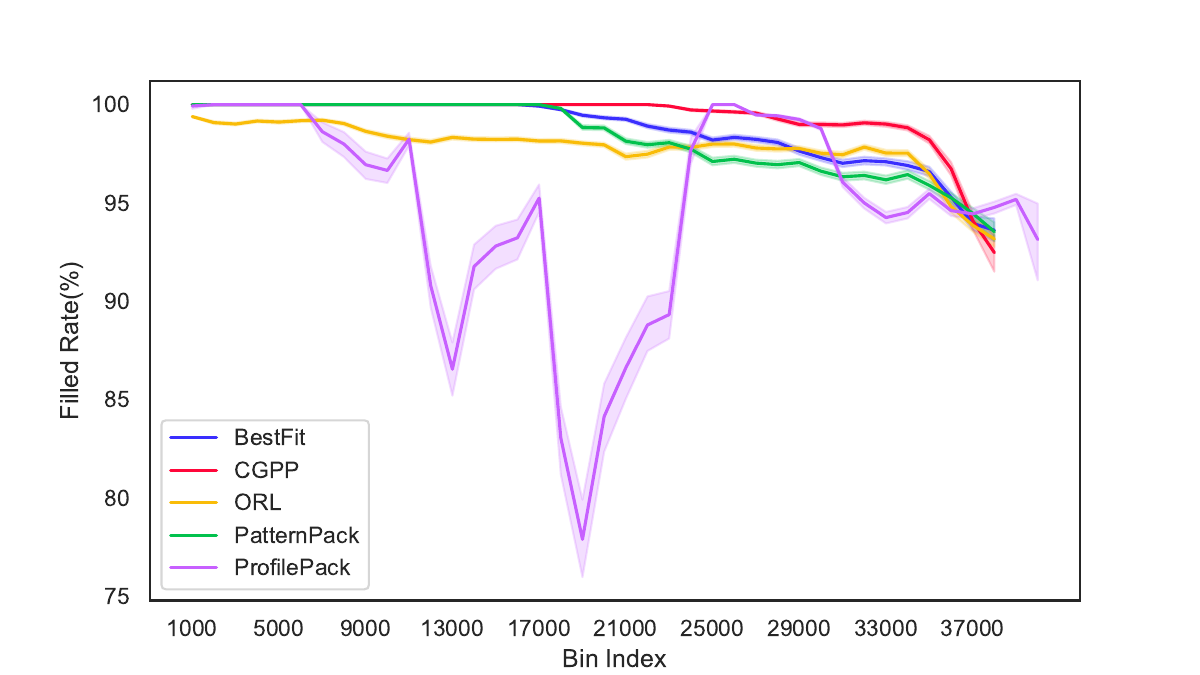}
\caption{Average Bin Filled Rate for Periodic Weibull}
\label{fig:solvis-periodic}
\end{figure}

Figure~\ref{fig:solvis-periodic} represents the filled rates of the Periodic Weibull dataset.
Most methods initially achieve a nearly 100\% filled rate, as the sequence is long enough to provide sufficient small items to fill the wasted space in the bins opened at early stage. Thanks to its forward-looking strategy in the form of patterns, \texttt{CGPP} maintained a high filled rate over the entire packing stages until the very end, indicating its success in adaptively identifying good patterns even as the distribution changes.
\texttt{ORL} still tended to sacrifice some waste space in order to achieve a more stable filled rate.
Both \texttt{BestFit} and \texttt{PatternPack} experienced a drop in filled rate in the middle stage of the packing. \texttt{PatternPack} has slightly worse fill rate than \texttt{BestFit}'s and significantly worse than \texttt{CGPP}'s, primarily due to imperfect prediction when the distribution changes.
\texttt{ProfilePack} exhibited instability, and its performance dropped significantly when the distribution changes.

The poor performance of \texttt{ProfilePack} and \texttt{PatternPack} in this case highlights the potential risk of poor prediction will misguide packing. \texttt{ProfilePack} utilised best fit descending, which is sensitive to distribution change. This resulted the most unstable behaviour in terms of filled rate.
% The poor performance of \texttt{ProfilePack} and \texttt{PatternPack} in this case highlights the potential risk of prediction-based algorithms in real-world data: the prediction can be far from realistic. 
% When the distribution changed, the prediction and associated patterns can be incorrect. 
% \texttt{ProfilePack} applied best-fit descending, which generates sub-optimal plans for the new distribution. 
% \texttt{PatternPack} allows all bins to periodically fallback,
\texttt{PatternPack} tracked the distribution for planning, and also allowed items being packed by best-fit heuristics regardless of the packing plan. These two strategy contributed to the robustness of the algorithm,
but the wrong plan can still be executed halfway, resulting in additional waste as illustrated. On the other hand, the pricing-based pattern identification remains robust under the new situation and therefore leads to the best performance among all other methods.
% The pricing-free pattern generation method, like best fit descending in 
% \texttt{ProfilePack} will generate patterns that . 

\subsection{Analysis of patterns and their reuse}

In this section, we analyse the detailed pattern quality and determine the extent to which the pattern contributes to achieving a good solution. We have selected an instance in Burke-4 (experiment 15) as a representative case for discussion. Similar behaviours can be observed from most other instances.

% We expected \texttt{CGPP} to recognise good patterns for guiding a superior solution. 

% The quality of a pattern is measured by its filled rate, which quantifies the ratio of the space occupied by the items in comparison to the bin's capacity. Moreover, patterns with higher quality were expected to occur more frequently in the final solution. 

We firstly provide an offline oracle solution with all information being known in advance. The bin patterns used in such an offline oracle solution might be regarded as high-quality patterns. We expect an algorithm that is able to recognise good patterns will tend to pack bins similarly to the oracle solution. That is, not only the high-quality patterns should be used more in the online solution, but also the pattern distribution should be close to the offline oracle.

\begin{figure}[t]
\centering
\includegraphics[width=\textwidth]{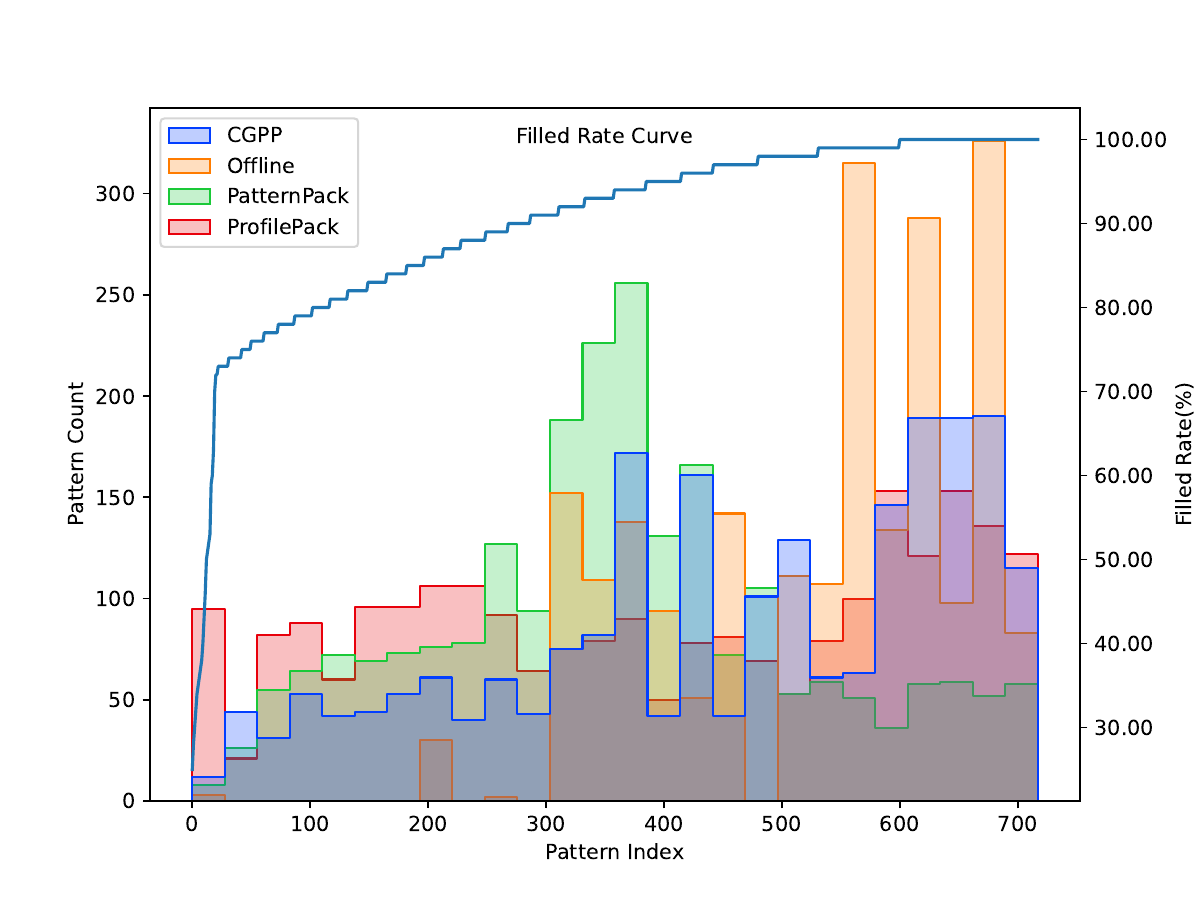}
\caption{Histogram of pattern quantity and their fill rates. Blue curve: fill rate of patterns measured by the second y-axis on the right.}
\label{fig:pattern}
\end{figure}

Figure~\ref{fig:pattern} represents the histogram of solution patterns. All patterns are sorted by their fill rates, and each pattern is assigned a unique index, where a larger index indicates a higher fill rate. The changes of fill rates across different pattern indices is represented by the blue curve in the figure (measured by the second y-axis on the right). The height of each histogram bar represents the quantity of a certain pattern used in the solution. The \texttt{Offline} histogram represents the pattern distribution for the offline oracle, where the patterns are considered in high-quality.
In comparison, \texttt{CGPP} achieved a histogram that closely resembles the offline solution, with more high-quality patterns being used and much higher overlap with oracle solution. This indicates that not only was \texttt{CGPP} able to identify good patterns, but it could also effectively reuse those patterns to reduce overall waste, resulting in improved bin usage.

For \texttt{PatternPack}, it also demonstrated the ability to reuse patterns. However, it favoured patterns at index 300-400, resulting in not only high frequency of sub-optimal patterns (90\%~95\% filled rate), but also low overlap with oracle patterns. In the case of \texttt{ProfilePack}, the patterns were more evenly distributed. It achieved better overlap at index 500-700 than \texttt{PatternPack} but it also used many patterns of low fill rates (e.g. pattern index 0-300), which rarely appeared in offline oracle. These low quality patterns resulted in poorer overall performance.

\subsection{Discussion on compared methods}
    % \subsubsection{Comparison with heuristics}
    % In most distribution in our experiment, \texttt{CGPP} consistently outperforms \texttt{BestFit}. Noting the fallback strategy of \texttt{CGPP} is \texttt{BestFit}, such advantage is mainly caused by the planning module. In other words, the online distribution learning performed well and the dualism pricing for identifying efficient patterns worked well. 
    In most datasets in our experiment, \texttt{CGPP} outperforms \texttt{BestFit} and other existing methods. The advantage can be mainly attributed to the dynamic pattern identification and associated planning, the reactive fallback strategy. Noting that with distribution of random variables being given, \texttt{CGPP} could achieve near-optimal solutions, significantly outperforming all other methods. Under unknown distributions, it achieved excellent performance in most cases when compared with other methods.

    % \subsubsection{Comparison with \texttt{ORL}} 
    % We also conducted a comparison between our method and 
    % \texttt{ORL} introduced in \citet{balaji_orlreinforcementlearning_2019}, the only reinforcement learning based online 1D bin packing work available as far as we know. In terms of solution quality, 
    Our proposed method significantly outperformed \texttt{ORL} approach. To our surprise, the \texttt{ORL} did not even outperform the \texttt{BestFit} strategy in most case. This could potentially be attributed to insufficient training time, considering the hardware limitations in our experiments compared to the original work. Furthermore, even when well-trained, the \texttt{ORL} exhibited weaker generalisation ability. As reported by \citet{balaji_orlreinforcementlearning_2019}, a RL model trained on PP/BW instances performed poorly on LW instances. In contrast, our method demonstrated strong generalisation capabilities, even with limited prior knowledge.
    
    % \subsubsection{Comparison with Plan Based Methods} 

    % The rationale behind comparing our method with \texttt{PatternPack} is the similar formulation of the problem: we both view the BPP as a sequential packing process following a pattern-based plan, where planning is executed online. 
    \texttt{PatternPack} is designed to solve problems with large discrete or continuous item types \citep{sym14071301}, which is not the case in this work. Our method demonstrated superior performance in terms of average bin usage compared to \texttt{PatternPack} and its fuzzy-enhanced version \texttt{FPP}, as our method utilised dualism pricing based Column Generation for planning, which typically yields better results with the online heuristic employed by \texttt{PatternPack}. Additionally, our uncertainty handling strategy can identify and eliminate planning errors caused by imperfect predictions. 
    % Regarding computational cost, although not directly compared through experiments, \texttt{PatternPack} performs better. \texttt{CGPP} heavily relies on column generation planning and re-planning during execution. This process still requires significant computational time as the problem volume increases due to the NP-Hard nature of the problem.

    \texttt{ProfilePack} theoretically proved that applying good prediction can lead to high-quality solutions. %is another work that theoretically proved that applying good prediction can lead to high-quality solutions. This work
    % is more focused on theoretical analysis. 
    In terms of implementation detail, the best-fit descending profile generation was not robust towards changing distributions. Also, the method lacks an uncertainty handling strategy. In some special cases (e.g., experiment 26), \texttt{ProfilePack} generates extremely poor results, indicating its major reliability issues. On the other hand, the success of \texttt{CGPP} in most cases also supports the theoretical result that good prediction can guide a better online packing strategy.
    
    % \subsection{Time Cost}
    % If the plan is perfect, the packing procedure has a time complexity of at most $O(N|P|\sum_{\mathbf{p}\in P} z_\mathbf{p})$, where $|P|$ is the total number of adopted patterns and $\sum_{\mathbf{p}\in P}z_\mathbf{p}$ is the total bin count. Considering that the total bin count is proportional to $N$, we can say that the time complexity is $O(N^2)$. When the plan is imperfect, the worst case will fallback to the best-fit heuristic, which also has an $O(N^2)$ time complexity.

    % The column generation planning involved in our method is NP-hard since it involves solving a linear integer programming problem for pattern generation and plan generation. 
    % % However, the computational complexity of the procedure is primarily dependent on the number of types, denoted as $T$. 
    % When the number of types is relatively small, the plan generation procedure can be considered constant. In our experience, it usually takes less than 10 seconds to complete the plan generation process for 100 types.

\section{Conclusion and future work}
In combinatorial optimisation, patterns are reusable building blocks of solutions that are more favourable than black-box solvers. However, we showed in this research that the values of patterns could change due to uncertainties related to objectives and constraints and most existing methods fail to exploit the inter-dependencies among decision variables incurred by uncertainties in constraints. We established a scheme to dynamically quantify the usefulness of different patterns based on the dualism of COP and use the information to guide the decision process. To handle the influence caused by both underestimation and overestimation, we introduced threshold-based methods to eliminate the inconsistency between plan and observation.
The test results on bin packing problem show significant performance advantage from the proposed method compared with the current the state-of-the-art methods. In future, we would investigate how the proposed framework generalises to other COPs with similar structures. %other novel ways to combine learning methods with patterns and dynamic pricing.
% A state monitor is used to monitor possible switches of distributions of the randomness, and bounding the uncertainty.


%%       \citet{<label>} ==> Jones et al. [21]
%%  ~\citep{<label>} ==> [21]
%%

 \section*{Acknowledgements}
This work was supported in part by the National Natural Science Foundation of China under Grant 72071116 and in part by the Ningbo Municipal Bureau of Science and Technology under Grant 2021Z173.

\small
\bibliographystyle{elsarticle-num-names} 
\bibliography{main.bib}

\begin{thebibliography}{44}
\expandafter\ifx\csname natexlab\endcsname\relax\def\natexlab#1{#1}\fi
\providecommand{\url}[1]{\texttt{#1}}
\providecommand{\href}[2]{#2}
\providecommand{\path}[1]{#1}
\providecommand{\DOIprefix}{doi:}
\providecommand{\ArXivprefix}{arXiv:}
\providecommand{\URLprefix}{URL: }
\providecommand{\Pubmedprefix}{pmid:}
\providecommand{\doi}[1]{\href{http://dx.doi.org/#1}{\path{#1}}}
\providecommand{\Pubmed}[1]{\href{pmid:#1}{\path{#1}}}
\providecommand{\bibinfo}[2]{#2}
\ifx\xfnm\relax \def\xfnm[#1]{\unskip,\space#1}\fi
%Type = Book
\bibitem[{Korte et~al.(2011)Korte, Vygen, Korte, and Vygen}]{korte2011combinatorial}
\bibinfo{author}{B.~H. Korte}, \bibinfo{author}{J.~Vygen}, \bibinfo{author}{B.~Korte}, \bibinfo{author}{J.~Vygen}, \bibinfo{title}{Combinatorial optimization}, volume~\bibinfo{volume}{1}, \bibinfo{publisher}{Springer}, \bibinfo{year}{2011}.
%Type = Article
\bibitem[{Bai et~al.(2023)Bai, Chen, Chen, Cui, Gong, He, Jiang, Jin, Jin, Kendall et~al.}]{bai-analytics-2023}
\bibinfo{author}{R.~Bai}, \bibinfo{author}{X.~Chen}, \bibinfo{author}{Z.-L. Chen}, \bibinfo{author}{T.~Cui}, \bibinfo{author}{S.~Gong}, \bibinfo{author}{W.~He}, \bibinfo{author}{X.~Jiang}, \bibinfo{author}{H.~Jin}, \bibinfo{author}{J.~Jin}, \bibinfo{author}{G.~Kendall}, et~al.,
\newblock \bibinfo{title}{Analytics and machine learning in vehicle routing research},
\newblock \bibinfo{journal}{International Journal of Production Research} \bibinfo{volume}{61} (\bibinfo{year}{2023}) \bibinfo{pages}{4--30}.
%Type = Article
\bibitem[{Gupta and Radovanovi{\'c}(2020)}]{gupta2020interior}
\bibinfo{author}{V.~Gupta}, \bibinfo{author}{A.~Radovanovi{\'c}},
\newblock \bibinfo{title}{Interior-point-based online stochastic bin packing},
\newblock \bibinfo{journal}{Operations Research} \bibinfo{volume}{68} (\bibinfo{year}{2020}) \bibinfo{pages}{1474--1492}.
%Type = Article
\bibitem[{Chen et~al.(2024)Chen, Bai, Qu, Dong, and Jin}]{chen2024deep}
\bibinfo{author}{X.~Chen}, \bibinfo{author}{R.~Bai}, \bibinfo{author}{R.~Qu}, \bibinfo{author}{J.~Dong}, \bibinfo{author}{Y.~Jin},
\newblock \bibinfo{title}{Deep reinforcement learning assisted genetic programming ensemble hyper-heuristics for dynamic scheduling of container port trucks},
\newblock \bibinfo{journal}{IEEE Transactions on Evolutionary Computation} \bibinfo{volume}{In press} (\bibinfo{year}{2024}). \DOIprefix\doi{10.1109/TEVC.2024.3381042}.
%Type = Inproceedings
\bibitem[{Lu et~al.(2020)Lu, Zhang, and Yang}]{Lu2020A}
\bibinfo{author}{H.~Lu}, \bibinfo{author}{X.~Zhang}, \bibinfo{author}{S.~Yang},
\newblock \bibinfo{title}{A learning-based iterative method for solving vehicle routing problems},
\newblock in: \bibinfo{booktitle}{International Conference on Learning Representations}, \bibinfo{year}{2020}.
%Type = Article
\bibitem[{Gomez et~al.(2014)Gomez, Alvarez, Jacobo-Berlles, and Mejail}]{GOMEZ20149}
\bibinfo{author}{L.~Gomez}, \bibinfo{author}{L.~Alvarez}, \bibinfo{author}{J.~Jacobo-Berlles}, \bibinfo{author}{M.~Mejail},
\newblock \bibinfo{title}{Special issue on computer vision applying pattern recognition techniques},
\newblock \bibinfo{journal}{Pattern Recognition} \bibinfo{volume}{47} (\bibinfo{year}{2014}) \bibinfo{pages}{9--11}.
%Type = Article
\bibitem[{Chen et~al.(2022)Chen, Wang, and Gao}]{CHEN2022108769}
\bibinfo{author}{X.~Chen}, \bibinfo{author}{B.~Wang}, \bibinfo{author}{Y.~Gao},
\newblock \bibinfo{title}{Symmetric binary tree based co-occurrence texture pattern mining for fine-grained plant leaf image retrieval},
\newblock \bibinfo{journal}{Pattern Recognition} \bibinfo{volume}{129} (\bibinfo{year}{2022}) \bibinfo{pages}{108769}.
%Type = Article
\bibitem[{Lin et~al.(2024)Lin, Jiang, and Zheng}]{LIN2024110143}
\bibinfo{author}{Z.~Lin}, \bibinfo{author}{X.~Jiang}, \bibinfo{author}{Z.~Zheng},
\newblock \bibinfo{title}{A coarse-to-fine pattern parser for mitigating the issue of drastic imbalance in pixel distribution},
\newblock \bibinfo{journal}{Pattern Recognition} \bibinfo{volume}{148} (\bibinfo{year}{2024}) \bibinfo{pages}{110143}.
%Type = Article
\bibitem[{Bai et~al.(2020)Bai, Hancock, Wilson, and Ho}]{BAI202046}
\bibinfo{author}{X.~Bai}, \bibinfo{author}{E.~R. Hancock}, \bibinfo{author}{R.~C. Wilson}, \bibinfo{author}{T.~K. Ho},
\newblock \bibinfo{title}{Special issue on recent advances in statistical, structural and syntactic pattern recognition},
\newblock \bibinfo{journal}{Pattern Recognition Letters} \bibinfo{volume}{131} (\bibinfo{year}{2020}) \bibinfo{pages}{46--48}.
%Type = Article
\bibitem[{Li et~al.(2021)Li, Zhao, Lee, Sassan, and Lin}]{LI2021107711}
\bibinfo{author}{Q.~Li}, \bibinfo{author}{L.~Zhao}, \bibinfo{author}{Y.-C. Lee}, \bibinfo{author}{A.~Sassan}, \bibinfo{author}{J.~Lin},
\newblock \bibinfo{title}{Cpm: A general feature dependency pattern mining framework for contrast multivariate time series},
\newblock \bibinfo{journal}{Pattern Recognition} \bibinfo{volume}{112} (\bibinfo{year}{2021}) \bibinfo{pages}{107711}.
%Type = Article
\bibitem[{Breitenbach et~al.(2023)Breitenbach, Wilkusz, Rasbach, and Jahnke}]{BREITENBACH2023109355}
\bibinfo{author}{T.~Breitenbach}, \bibinfo{author}{B.~Wilkusz}, \bibinfo{author}{L.~Rasbach}, \bibinfo{author}{P.~Jahnke},
\newblock \bibinfo{title}{On a method for detecting periods and repeating patterns in time series data with autocorrelation and function approximation},
\newblock \bibinfo{journal}{Pattern Recognition} \bibinfo{volume}{138} (\bibinfo{year}{2023}) \bibinfo{pages}{109355}.
%Type = Article
\bibitem[{Martello and Toth(1990)}]{martello_lowerboundsreduction_1990}
\bibinfo{author}{S.~Martello}, \bibinfo{author}{P.~Toth},
\newblock \bibinfo{title}{Lower bounds and reduction procedures for the bin packing problem},
\newblock \bibinfo{journal}{Discrete Applied Mathematics} \bibinfo{volume}{28} (\bibinfo{year}{1990}) \bibinfo{pages}{59--70}.
%Type = Article
\bibitem[{Bai et~al.(2012)Bai, Blazewicz, Burke, Kendall, and McCollum}]{bai-simulated-2012}
\bibinfo{author}{R.~Bai}, \bibinfo{author}{J.~Blazewicz}, \bibinfo{author}{E.~Burke}, \bibinfo{author}{G.~Kendall}, \bibinfo{author}{B.~McCollum},
\newblock \bibinfo{title}{A simulated annealing hyper-heuristic methodology for flexible decision support},
\newblock \bibinfo{journal}{4OR} \bibinfo{volume}{10} (\bibinfo{year}{2012}) \bibinfo{pages}{43--66}.
%Type = Book
\bibitem[{Scheithauer(2018)}]{Scheithauer2018}
\bibinfo{author}{G.~Scheithauer}, \bibinfo{title}{Introduction to Cutting and Packing Optimization: Problems, Modeling Approaches, Solution Methods}, \bibinfo{publisher}{Springer}, \bibinfo{year}{2018}.
%Type = Book
\bibitem[{Bonnans(2019)}]{bonnans_convex_2019}
\bibinfo{author}{J.~F. Bonnans}, \bibinfo{title}{Convex and Stochastic Optimization}, \bibinfo{publisher}{Springer}, \bibinfo{year}{2019}.
%Type = Article
\bibitem[{Kuosmanen and Zhou(2021)}]{KUOSMANEN2021666}
\bibinfo{author}{T.~Kuosmanen}, \bibinfo{author}{X.~Zhou},
\newblock \bibinfo{title}{Shadow prices and marginal abatement costs: Convex quantile regression approach},
\newblock \bibinfo{journal}{European Journal of Operational Research} \bibinfo{volume}{289} (\bibinfo{year}{2021}) \bibinfo{pages}{666--675}.
%Type = Inproceedings
\bibitem[{Li et~al.(2021)Li, Xia, Yan, Sun, Zhao, and Liu}]{Jinpeng-Stylized-2021}
\bibinfo{author}{J.~Li}, \bibinfo{author}{Y.~Xia}, \bibinfo{author}{R.~Yan}, \bibinfo{author}{H.~Sun}, \bibinfo{author}{D.~Zhao}, \bibinfo{author}{T.-Y. Liu},
\newblock \bibinfo{title}{Stylized dialogue generation with multi-pass dual learning},
\newblock in: \bibinfo{editor}{M.~Ranzato}, \bibinfo{editor}{A.~Beygelzimer}, \bibinfo{editor}{Y.~Dauphin}, \bibinfo{editor}{P.~Liang}, \bibinfo{editor}{J.~W. Vaughan} (Eds.), \bibinfo{booktitle}{Advances in Neural Information Processing Systems}, volume~\bibinfo{volume}{34}, \bibinfo{publisher}{Curran Associates, Inc.}, \bibinfo{year}{2021}, pp. \bibinfo{pages}{28470--28481}.
%Type = Article
\bibitem[{Wang et~al.(2019)Wang, Xia, He, Tian, Qin, Zhai, and Liu}]{wang2019MultiDualLearn}
\bibinfo{author}{Y.~Wang}, \bibinfo{author}{Y.~Xia}, \bibinfo{author}{T.~He}, \bibinfo{author}{F.~Tian}, \bibinfo{author}{T.~Qin}, \bibinfo{author}{C.~Zhai}, \bibinfo{author}{T.-Y. Liu},
\newblock \bibinfo{title}{Multi-{{Agent Dual Learning}}},
\newblock \bibinfo{journal}{Proceedings of the International Conference on Learning Representations (ICLR) 2019}  (\bibinfo{year}{2019}).
%Type = Article
\bibitem[{Delorme et~al.(2016)Delorme, Iori, and Martello}]{DELORME20161}
\bibinfo{author}{M.~Delorme}, \bibinfo{author}{M.~Iori}, \bibinfo{author}{S.~Martello},
\newblock \bibinfo{title}{Bin packing and cutting stock problems: Mathematical models and exact algorithms},
\newblock \bibinfo{journal}{European Journal of Operational Research} \bibinfo{volume}{255} (\bibinfo{year}{2016}) \bibinfo{pages}{1--20}.
%Type = Inproceedings
\bibitem[{Wang et~al.(2010)Wang, Luh, Gribik, Zhang, and Peng}]{wang_subgradientsimplexbasedcutting_2010}
\bibinfo{author}{C.~Wang}, \bibinfo{author}{P.~B. Luh}, \bibinfo{author}{P.~Gribik}, \bibinfo{author}{L.~Zhang}, \bibinfo{author}{T.~Peng},
\newblock \bibinfo{title}{The subgradient-simplex based cutting plane method for convex hull pricing},
\newblock in: \bibinfo{booktitle}{{{IEEE PES General Meeting}}}, \bibinfo{year}{2010}, pp. \bibinfo{pages}{1--8}.
%Type = Article
\bibitem[{Xue et~al.(2021)Xue, Bai, Qu, and Aickelin}]{xue_hybridpricingcutting_2021}
\bibinfo{author}{N.~Xue}, \bibinfo{author}{R.~Bai}, \bibinfo{author}{R.~Qu}, \bibinfo{author}{U.~Aickelin},
\newblock \bibinfo{title}{A hybrid pricing and cutting approach for the multi-shift full truckload vehicle routing problem},
\newblock \bibinfo{journal}{European Journal of Operational Research} \bibinfo{volume}{292} (\bibinfo{year}{2021}) \bibinfo{pages}{500--514}.
%Type = Article
\bibitem[{Sheng et~al.(2022)Sheng, Hu, Zhou, Zhu, Jin, Wang, and Wang}]{Junjie-learning-2021}
\bibinfo{author}{J.~Sheng}, \bibinfo{author}{Y.~Hu}, \bibinfo{author}{W.~Zhou}, \bibinfo{author}{L.~Zhu}, \bibinfo{author}{B.~Jin}, \bibinfo{author}{J.~Wang}, \bibinfo{author}{X.~Wang},
\newblock \bibinfo{title}{Learning to schedule multi-numa virtual machines via reinforcement learning},
\newblock \bibinfo{journal}{Pattern Recognition} \bibinfo{volume}{121} (\bibinfo{year}{2022}) \bibinfo{pages}{108254}.
%Type = Article
\bibitem[{Abdalkareem et~al.(2021)Abdalkareem, Amir, {Al-Betar}, Ekhan, and Hammouri}]{abdalkareem_healthcareschedulingoptimization_2021c}
\bibinfo{author}{Z.~A. Abdalkareem}, \bibinfo{author}{A.~Amir}, \bibinfo{author}{M.~A. {Al-Betar}}, \bibinfo{author}{P.~Ekhan}, \bibinfo{author}{A.~I. Hammouri},
\newblock \bibinfo{title}{Healthcare scheduling in optimization context: A review},
\newblock \bibinfo{journal}{Health and Technology} \bibinfo{volume}{11} (\bibinfo{year}{2021}) \bibinfo{pages}{445--469}.
%Type = Article
\bibitem[{Ali et~al.(2022)Ali, Ramos, Carravilla, and Oliveira}]{ALI2022108122}
\bibinfo{author}{S.~Ali}, \bibinfo{author}{A.~G. Ramos}, \bibinfo{author}{M.~A. Carravilla}, \bibinfo{author}{J.~F. Oliveira},
\newblock \bibinfo{title}{On-line three-dimensional packing problems: A review of off-line and on-line solution approaches},
\newblock \bibinfo{journal}{Computers \& Industrial Engineering} \bibinfo{volume}{168} (\bibinfo{year}{2022}) \bibinfo{pages}{108122}.
%Type = Incollection
\bibitem[{Coffman et~al.(2013)Coffman, Csirik, Galambos, Martello, and Vigo}]{coffman_binpackingapproximation_2013}
\bibinfo{author}{E.~G. Coffman}, \bibinfo{author}{J.~Csirik}, \bibinfo{author}{G.~Galambos}, \bibinfo{author}{S.~Martello}, \bibinfo{author}{D.~Vigo},
\newblock \bibinfo{title}{Bin {{Packing Approximation Algorithms}}: {{Survey}} and {{Classification}}},
\newblock in: \bibinfo{editor}{P.~M. Pardalos}, \bibinfo{editor}{D.-Z. Du}, \bibinfo{editor}{R.~L. Graham} (Eds.), \bibinfo{booktitle}{Handbook of {{Combinatorial Optimization}}}, \bibinfo{publisher}{{Springer New York}}, \bibinfo{address}{{New York, NY}}, \bibinfo{year}{2013}, pp. \bibinfo{pages}{455--531}.
%Type = Article
\bibitem[{Zhang et~al.(2020)Zhang, Liu, Zhang, and Qin}]{Zhang2020ColumnGA}
\bibinfo{author}{Q.~Zhang}, \bibinfo{author}{S.~Liu}, \bibinfo{author}{R.~Zhang}, \bibinfo{author}{S.~Qin},
\newblock \bibinfo{title}{Column generation algorithms for mother plate design in steel plants},
\newblock \bibinfo{journal}{OR Spectrum} \bibinfo{volume}{43} (\bibinfo{year}{2020}) \bibinfo{pages}{127--153}.
%Type = Article
\bibitem[{Dell’Amico et~al.(2020)Dell’Amico, Furini, and Iori}]{DELLAMICO2020104825}
\bibinfo{author}{M.~Dell’Amico}, \bibinfo{author}{F.~Furini}, \bibinfo{author}{M.~Iori},
\newblock \bibinfo{title}{A branch-and-price algorithm for the temporal bin packing problem},
\newblock \bibinfo{journal}{Computers \& Operations Research} \bibinfo{volume}{114} (\bibinfo{year}{2020}) \bibinfo{pages}{104825}.
%Type = Inproceedings
\bibitem[{Burke et~al.(2010)Burke, Hyde, and Kendall}]{5586388}
\bibinfo{author}{E.~K. Burke}, \bibinfo{author}{M.~R. Hyde}, \bibinfo{author}{G.~Kendall},
\newblock \bibinfo{title}{Providing a memory mechanism to enhance the evolutionary design of heuristics},
\newblock in: \bibinfo{booktitle}{IEEE Congress on Evolutionary Computation}, \bibinfo{year}{2010}, pp. \bibinfo{pages}{1--8}.
%Type = Article
\bibitem[{{L{\'o}pez-Camacho} et~al.(2014){L{\'o}pez-Camacho}, {Terashima-Marin}, Ross, and Ochoa}]{lopez-camacho_unifiedhyperheuristicframework_2014}
\bibinfo{author}{E.~{L{\'o}pez-Camacho}}, \bibinfo{author}{H.~{Terashima-Marin}}, \bibinfo{author}{P.~Ross}, \bibinfo{author}{G.~Ochoa},
\newblock \bibinfo{title}{A unified hyper-heuristic framework for solving bin packing problems},
\newblock \bibinfo{journal}{Expert Systems with Applications} \bibinfo{volume}{41} (\bibinfo{year}{2014}) \bibinfo{pages}{6876--6889}.
%Type = Article
\bibitem[{Liu et~al.(2021)Liu, Cheng, Tian, Wang, Leng, Zhao, Zhang, and Wei}]{LIU2021107175}
\bibinfo{author}{Q.~Liu}, \bibinfo{author}{H.~Cheng}, \bibinfo{author}{T.~Tian}, \bibinfo{author}{Y.~Wang}, \bibinfo{author}{J.~Leng}, \bibinfo{author}{R.~Zhao}, \bibinfo{author}{H.~Zhang}, \bibinfo{author}{L.~Wei},
\newblock \bibinfo{title}{Algorithms for the variable-sized bin packing problem with time windows},
\newblock \bibinfo{journal}{Computers \& Industrial Engineering} \bibinfo{volume}{155} (\bibinfo{year}{2021}) \bibinfo{pages}{107175}.
%Type = Article
\bibitem[{Yan et~al.(2022)Yan, Weng, Huang, Li, Zhou, Su, and Zhu}]{Dong-deep-2021}
\bibinfo{author}{D.~Yan}, \bibinfo{author}{J.~Weng}, \bibinfo{author}{S.~Huang}, \bibinfo{author}{C.~Li}, \bibinfo{author}{Y.~Zhou}, \bibinfo{author}{H.~Su}, \bibinfo{author}{J.~Zhu},
\newblock \bibinfo{title}{Deep reinforcement learning with credit assignment for combinatorial optimization},
\newblock \bibinfo{journal}{Pattern Recognition} \bibinfo{volume}{124} (\bibinfo{year}{2022}) \bibinfo{pages}{108466}.
%Type = Inproceedings
\bibitem[{Angelopoulos et~al.(2022)Angelopoulos, Kamali, and Shadkami}]{angelopoulos2022online}
\bibinfo{author}{S.~Angelopoulos}, \bibinfo{author}{S.~Kamali}, \bibinfo{author}{K.~Shadkami},
\newblock \bibinfo{title}{Online bin packing with predictions},
\newblock in: \bibinfo{booktitle}{Thirty-First International Joint Conference on Artificial Intelligence $\{$IJCAI-22$\}$}, volume~\bibinfo{volume}{36}, \bibinfo{organization}{International Joint Conferences on Artificial Intelligence Organization}, \bibinfo{year}{2022}, pp. \bibinfo{pages}{4574--4580}.
%Type = Article
\bibitem[{Lin et~al.(2024)Lin, Li, Cui, Jin, Bai, Qu, and Garibaldi}]{lin2024PatteAlgorFuzzy}
\bibinfo{author}{B.~Lin}, \bibinfo{author}{J.~Li}, \bibinfo{author}{T.~Cui}, \bibinfo{author}{H.~Jin}, \bibinfo{author}{R.~Bai}, \bibinfo{author}{R.~Qu}, \bibinfo{author}{J.~Garibaldi},
\newblock \bibinfo{title}{A pattern-based algorithm with fuzzy logic bin selector for online bin packing problem},
\newblock \bibinfo{journal}{Expert Systems with Applications} \bibinfo{volume}{249} (\bibinfo{year}{2024}) \bibinfo{pages}{123515}.
%Type = Article
\bibitem[{Bengio et~al.(2021)Bengio, Lodi, and Prouvost}]{BENGIO2021405}
\bibinfo{author}{Y.~Bengio}, \bibinfo{author}{A.~Lodi}, \bibinfo{author}{A.~Prouvost},
\newblock \bibinfo{title}{Machine learning for combinatorial optimization: A methodological tour d’horizon},
\newblock \bibinfo{journal}{European Journal of Operational Research} \bibinfo{volume}{290} (\bibinfo{year}{2021}) \bibinfo{pages}{405--421}.
%Type = Article
\bibitem[{Dai et~al.(2018)Dai, Khalil, Zhang, Dilkina, and Song}]{dai_learningcombinatorialoptimization_2018}
\bibinfo{author}{H.~Dai}, \bibinfo{author}{E.~B. Khalil}, \bibinfo{author}{Y.~Zhang}, \bibinfo{author}{B.~Dilkina}, \bibinfo{author}{L.~Song},
\newblock \bibinfo{title}{Learning {{Combinatorial Optimization Algorithms}} over {{Graphs}}},
\newblock \bibinfo{journal}{arXiv:1704.01665 [cs, stat]}  (\bibinfo{year}{2018}). \href{http://arxiv.org/abs/1704.01665}{{\tt arXiv:1704.01665}}.
%Type = Article
\bibitem[{Zhang et~al.(2021)Zhang, Zi, and Ge}]{zhang_attend2packbinpacking_2021}
\bibinfo{author}{J.~Zhang}, \bibinfo{author}{B.~Zi}, \bibinfo{author}{X.~Ge},
\newblock \bibinfo{title}{{{Attend2Pack}}: {{Bin Packing}} through {{Deep Reinforcement Learning}} with {{Attention}}}  (\bibinfo{year}{2021}). \href{http://arxiv.org/abs/2107.04333}{{\tt arXiv:2107.04333}}.
%Type = Article
\bibitem[{Hubbs et~al.(2020)Hubbs, Perez, Sarwar, Sahinidis, Grossmann, and Wassick}]{hubbs_orgymreinforcementlearning_2020}
\bibinfo{author}{C.~D. Hubbs}, \bibinfo{author}{H.~D. Perez}, \bibinfo{author}{O.~Sarwar}, \bibinfo{author}{N.~V. Sahinidis}, \bibinfo{author}{I.~E. Grossmann}, \bibinfo{author}{J.~M. Wassick},
\newblock \bibinfo{title}{Or-gym: {A} reinforcement learning library for operations research problem},
\newblock \bibinfo{journal}{CoRR} \bibinfo{volume}{abs/2008.06319} (\bibinfo{year}{2020}). \href{http://arxiv.org/abs/2008.06319}{{\tt arXiv:2008.06319}}.
%Type = Article
\bibitem[{Balaji et~al.(2019)Balaji, Bell{-}Masterson, Bilgin, Damianou, Garcia, Jain, Luo, Maggiar, Narayanaswamy, and Ye}]{balaji_orlreinforcementlearning_2019}
\bibinfo{author}{B.~Balaji}, \bibinfo{author}{J.~Bell{-}Masterson}, \bibinfo{author}{E.~Bilgin}, \bibinfo{author}{A.~C. Damianou}, \bibinfo{author}{P.~M. Garcia}, \bibinfo{author}{A.~Jain}, \bibinfo{author}{R.~Luo}, \bibinfo{author}{A.~Maggiar}, \bibinfo{author}{B.~Narayanaswamy}, \bibinfo{author}{C.~Ye},
\newblock \bibinfo{title}{{ORL:} reinforcement learning benchmarks for online stochastic optimization problems},
\newblock \bibinfo{journal}{CoRR} \bibinfo{volume}{abs/1911.10641} (\bibinfo{year}{2019}).
%Type = Inproceedings
\bibitem[{Zhao et~al.(2021)Zhao, She, Zhu, Yang, and Xu}]{zhao_online3dbin_2021}
\bibinfo{author}{H.~Zhao}, \bibinfo{author}{Q.~She}, \bibinfo{author}{C.~Zhu}, \bibinfo{author}{Y.~Yang}, \bibinfo{author}{K.~Xu},
\newblock \bibinfo{title}{Online 3d bin packing with constrained deep reinforcement learning},
\newblock in: \bibinfo{booktitle}{Thirty-Fifth {AAAI} Conference on Artificial Intelligence, {AAAI}}, \bibinfo{year}{2021}, pp. \bibinfo{pages}{741--749}.
%Type = Inproceedings
\bibitem[{Zhao and Xu(2022)}]{zhao2022learning}
\bibinfo{author}{H.~Zhao}, \bibinfo{author}{K.~Xu},
\newblock \bibinfo{title}{Learning efficient online 3d bin packing on packing configuration trees},
\newblock in: \bibinfo{booktitle}{International Conference on Learning Representations}, \bibinfo{year}{2022}.
%Type = Inproceedings
\bibitem[{Gao et~al.(2020)Gao, Wang, and Shen}]{9209730}
\bibinfo{author}{J.~Gao}, \bibinfo{author}{H.~Wang}, \bibinfo{author}{H.~Shen},
\newblock \bibinfo{title}{Machine learning based workload prediction in cloud computing},
\newblock in: \bibinfo{booktitle}{2020 29th International Conference on Computer Communications and Networks (ICCCN)}, \bibinfo{year}{2020}, pp. \bibinfo{pages}{1--9}.
%Type = Article
\bibitem[{Lin et~al.(2022)Lin, Li, Bai, Qu, Cui, and Jin}]{sym14071301}
\bibinfo{author}{B.~Lin}, \bibinfo{author}{J.~Li}, \bibinfo{author}{R.~Bai}, \bibinfo{author}{R.~Qu}, \bibinfo{author}{T.~Cui}, \bibinfo{author}{H.~Jin},
\newblock \bibinfo{title}{Identify patterns in online bin packing problem: An adaptive pattern-based algorithm},
\newblock \bibinfo{journal}{Symmetry} \bibinfo{volume}{14} (\bibinfo{year}{2022}).
%Type = Inproceedings
\bibitem[{Burke et~al.(2010)Burke, Hyde, and Kendall}]{burke_providingmemorymechanism_2010}
\bibinfo{author}{E.~K. Burke}, \bibinfo{author}{M.~R. Hyde}, \bibinfo{author}{G.~Kendall},
\newblock \bibinfo{title}{Providing a memory mechanism to enhance the evolutionary design of heuristics},
\newblock in: \bibinfo{booktitle}{{{IEEE Congress}} on {{Evolutionary Computation}}}, \bibinfo{year}{2010}, pp. \bibinfo{pages}{1--8}.
%Type = Inproceedings
\bibitem[{Casti{\~n}eiras et~al.(2012)Casti{\~n}eiras, De~Cauwer, and O'Sullivan}]{casti2012WeibuBenchBin}
\bibinfo{author}{I.~Casti{\~n}eiras}, \bibinfo{author}{M.~De~Cauwer}, \bibinfo{author}{B.~O'Sullivan},
\newblock \bibinfo{title}{Weibull-{{Based Benchmarks}} for {{Bin Packing}}},
\newblock in: \bibinfo{editor}{M.~Milano} (Ed.), \bibinfo{booktitle}{Principles and {{Practice}} of {{Constraint Programming}}}, Lecture {{Notes}} in {{Computer Science}}, \bibinfo{publisher}{Springer}, \bibinfo{address}{Berlin, Heidelberg}, \bibinfo{year}{2012}, pp. \bibinfo{pages}{207--222}.

\end{thebibliography}





% \end{thebibliography}

\end{document}